\documentclass[12pt]{article}
\author{Fan Kong, Keyan Song, Pu Zhang}
\title{Decomposition of \text{Torsion Pair}s on Module Categories}
\date{}
\usepackage{amsfonts}
\usepackage{graphicx}
\usepackage{amscd}
\usepackage[all]{xy}
\usepackage{amsmath,  cchead,  amsfonts,  amssymb,  amsthm,  mathrsfs,  xypic}
\usepackage{leftidx}
\usepackage[center]{titlesec}
\usepackage{ulem}
\usepackage{chemarrow}
 \xyoption{curve}
\newtheorem{thm}{Theorem}[section]
\newtheorem{cor}[thm]{Corollary}

\newtheorem{lem}[thm]{Lemma}
\newtheorem{exm}[thm]{Example}

\newtheorem{prop}[thm]{Proposition}
\newtheorem{rem}[thm]{Remark}
\newtheorem{defn}[thm]{Definition}
\numberwithin{equation}{section}

\newcommand{\ra}{\rightarrow}

\newcommand{\A}{\Lambda}

\newcommand{\add}{\operatorname{add}}
\newcommand{\Ima}{\operatorname{Im}}
\newcommand{\Cok}{\operatorname{Coker}}
\newcommand{\Ker}{\operatorname{Ker}}
\newcommand{\Hom}{\operatorname{Hom}}
\newcommand{\Ext}{\operatorname{Ext}}

\begin{document}

\maketitle

\begin{flushleft}
{\bf Abstract:} In this article, we generalize the concept of \text{torsion pair}s and study its structure. As a trial of obtaining all \text{torsion pair}s, we decompose \text{torsion pair}s by projective modules and injective modules. Then we calculate \text{torsion pair}s on the algebra $KA_n$ and tub categories. At last we try to find all torsion pairs on the module categories of finite dimensional hereditary algebras.

\vskip5pt

{\bf Key words:}  n-\text{torsion pair}, n-\text{torsion pair} seires, 1-type part partition, 2-type part partition, $\Ext$-projective, $\Ext$-injective.
\end{flushleft}

\vskip10pt
\section{Introduction}
The concept of torsion pair  on abelian category was introduced by Dickson in 1966 \cite{10}. From that time on, \text{torsion pair} has been always  a useful tool for studying the structure of module categories. However,  it seems there is no useful way to find all torsion pairs of a given algebra, although indeed there are some ways to construct torsion pairs among which the most well known is the tilting theory. As a trial, we try to give a way to obtain all torsion pairs of hereditary algebras in this article.  This topic is also discussed by Assem and Kerner in  \cite{1} where their most interest is to classify and characterize the torsion pairs by partial tilting modules.

   In  section 2, we study the general theory where  we introduce $n$-\text{torsion pair} and $n$-\text{torsion pair} series as the generalization of classic \text{torsion pair} and study its structure.  We can see that these two generalizations are essentially the same. In the rest of the paper we would know it is necessary and natural to put forward this conception  for studying the structure of torsion pairs. The main skill in this section is from \cite{8} and \cite{12} where they study HN-filtration for some categories.  There are really a lot of examples to illustrate the necessity to study this finer structure of module categories. For example perpendicular category is obtained by a 2-\text{torsion pair} series, and the structure of partial tilting  modules can be considered in this way. And  HN-filtration can be seen as a generalized $n-\text{torsion pair}$.

In \cite{1}, Assem and Kerner show a relation between some particular partial tilting moules and \text{torsion pair}s. In section 3, we adopt their ways by restricting to projective modules and injective modules to try to decompose all \text{torsion pair}s. And this is also an application of theories developed in section 2. We give a method for how to decompose a classic \text{torsion pair} to $n$-\text{torsion pair}s, and we give a  one to one  correspondence between all the \text{torsion pair}s and some sepcial $n$-\text{torsion pair} on the module category of any artin algebra.

In section 4, we apply the theory in section 3  to path algebras. As a application, we give all the \text{torsion pair}s on path algebra $KA_n$ and tube categories. Some of the results also have been shown in \cite{2} and  \cite{11}. But we think our results will be much more clear in some aspects.

The section 5 is devoted to obtain all torsion pairs of hereditary algebras which is our purpose. We define an operation called the translation of torsion pairs.  Combining this with the operation developed in section 3 and 4, the issue of obtaining all torsion pairs comes down to find all torsion pairs on regular component. For tame hereditary algebras, this problem is equivalent to calculate all torsion pairs on the tube categories in section 4.

We should admit that our  way of obtaining all torsion pairs is not very satisfactory since it is mixed with $\operatorname{DTr}$-translation and the extension between different parts of $n$-torsion pairs.

If there is no special instruction, all modules are left finitely generated modules. For an artin algebra $\A$, we denote by $\Lambda\text{-mod}$ the category of all left finitely generated $\Lambda$-modules. Subcategories are always assumed to be closed under isomorphism.

\section{$n-\text{torsion pair}$ and $n-\text{torsion pair}$ series}

In this section, we assume that $\Lambda$ is an artin algebra and  $\mathcal{C}$ is an extension-closed full subcategory of $\Lambda\text{-mod}$. If ${\mathcal{C}}_1,{\mathcal{C}}_2,\cdots,{\mathcal{C}}_n$ are full subcategories of  $\Lambda\text{-mod}$,  then we denote the minimal full extension-closed subcategory containing ${\mathcal{C}}_1,{\mathcal{C}}_2,\cdots,{\mathcal{C}}_n$ by $\langle {\mathcal{C}}_1,{\mathcal{C}}_2,\cdots,{\mathcal{C}}_n\rangle $.  If $\mathcal{D}$ is a subcategory of $\A$-mod, then we denote the set $\{M \mid \operatorname{Hom}(M,N) = 0, \forall N\in \mathcal{D}\} $  by $\leftidx{^\bot} \!\mathcal{D}$, the set$\{N \mid \operatorname{Hom}(M,N) = 0, \forall M\in \mathcal{D}\}$ by $\mathcal{D}^\bot$. \\The following definition is well known but different from that in \cite{13}.

 \begin{defn} A pair $(\mathcal{T},\mathcal{F})$ of full subcategories of $\mathcal{C}$ is called a \text{torsion pair} on $\mathcal{C}$  if the following conditions are satisfied:\\
$\left( 1 \right)$ $\operatorname{Hom}(X,Y) = 0$ for all $X\in\mathcal{T}$, $Y\in\mathcal{F}$.\\
$\left( 2 \right)$  $\forall X\in\mathcal{C}$, there exists an exact sequence on  $\Lambda\text{-mod}$:
$$ 0\longrightarrow X_{\mathcal{T}} \longrightarrow X \longrightarrow X_{\mathcal{F}}\longrightarrow 0$$
such that $X_{\mathcal{T}}\in \mathcal{T} $ and $X_{\mathcal{F}}\in {\mathcal{F}}.$
 \end{defn}

 \begin{rem}
 Let $(\mathcal{T},\mathcal{F})$ be a \text{torsion pair} on $\mathcal{C}$. Then
 $\mathcal{T} = \leftidx{^\bot} \!\mathcal{F} \bigcap \mathcal{C}$;  $\mathcal{F} = {\mathcal{T}}^\bot \bigcap \mathcal{C}$;  $\mathcal{T}$ and $\mathcal{F}$ are closed under extensions.
\end{rem}

Now we give the following definition which is a generalization of the above.

 \begin{defn}  an $n$-tuple $({\mathcal{C}}_1,{\mathcal{C}}_2,\cdots,{\mathcal{C}}_{n+1})$ of full extension-closed subcategories of $\mathcal{C}$ is called an $n-\text{torsion pair}$ if the following conditions are satisfied.\\
$\left( 1 \right)$  ${\mathcal{C}}_i  =  \mathcal{C} \bigcap  { \langle  {\mathcal{C}}_1,\cdots,{\mathcal{C}}_{i-1} \rangle }^{\bot}\bigcap  \leftidx{^{\bot}}  { \langle  {\mathcal{C}}_{i+1},\cdots,{\mathcal{C}}_{n+1} \rangle }$ for $i = 1,2,\cdots,n+1$.\\
$\left( 2 \right)$ $(\langle  {\mathcal{C}}_1,\cdots,{\mathcal{C}}_{i} \rangle ,\langle  {\mathcal{C}}_{i+1},\cdots,{\mathcal{C}}_{n+1} \rangle )$ is a $\text{torsion pair}$ on $\mathcal{C}$ for  $i = 1,2,\cdots,n+1$.\\
Moreover, if the first condition does not satisfy, we call $({\mathcal{C}}_1,{\mathcal{C}}_2,\cdots,{\mathcal{C}}_{n+1})$ a
defect $n-\text{torsion pair}$ on $\mathcal{C}$.
 \end{defn}

 The following lemma is obvious.
\begin{lem} Let ${\mathcal{C}}_1,{\mathcal{C}}_2,{\mathcal{C}}_3$ be $3$ full subcategories of $\Lambda\text{-mod}$. Then\\
$\left( 1 \right)$  $\langle {\mathcal{C}}_1,{\mathcal{C}}_2,{\mathcal{C}}_3\rangle   =  \langle \langle {\mathcal{C}}_1,{\mathcal{C}}_2\rangle ,{\mathcal{C}}_3\rangle   =  \langle {\mathcal{C}}_1,\langle {\mathcal{C}}_2,{\mathcal{C}}_3\rangle \rangle $.\\
$\left( 2 \right)$  $ \leftidx{^{\bot}}  { \langle  {\mathcal{C}}_1,{\mathcal{C}}_2 \rangle }  =   \leftidx{^{\bot}} { {\mathcal{C}}_1}\bigcap  \leftidx{^{\bot}} { {\mathcal{C}}_2} $,
${ \langle  {\mathcal{C}}_1,{\mathcal{C}}_2 \rangle }^{\bot}  =  {{\mathcal{C}}_1}^{\bot} \bigcap {{\mathcal{C}}_2}^{\bot}$.\\
$\left( 3 \right)$  $ \leftidx{^{\bot}}  { \langle  {\mathcal{C}}_1\rangle }  =   \leftidx{^{\bot}}  {  {\mathcal{C}}_1}$,${ \langle  {\mathcal{C}}_1 \rangle }^{\bot}  =  {  {\mathcal{C}}_1 }^{\bot}$.
\end{lem}

\begin{prop} Let
$({\mathcal{C}}_1,{\mathcal{C}}_2,\cdots,{\mathcal{C}}_{n+1})$ be an $n-\text{torsion pair}$ on $\mathcal{C}$. If $(\tilde{\mathcal{C}}_1,
\tilde{\mathcal{C}}_2,\cdots,\tilde{\mathcal{C}}_{k+1})$ be a $k-\text{torsion pair}$
on ${\mathcal{C}}_i$ for some $i$. Then
$({\mathcal{C}}_1,\cdots,{\mathcal{C}}_{i-1},
 \tilde{\mathcal{C}}_1,
\cdots,\tilde{\mathcal{C}}_{k+1},{\mathcal{C}}_{i+1},\cdots,{\mathcal{C}}_{n+1})$
is a $(n+k)-\text{torsion pair}$ on $\mathcal{C}$.
\end{prop}
\noindent{\bf Proof:} Step 1. If ${\mathcal{C}}_s \in {\{{\mathcal{C}}_1,{\mathcal{C}}_2,\cdots,{\mathcal{C}}_{i-1}\}}$,then

  $\mathcal{C}\bigcap  { \langle  {\mathcal{C}}_1,\cdots,{\mathcal{C}}_{s-1} \rangle }^{\bot}\bigcap  \leftidx{^{\bot}}  { \langle  {\mathcal{C}}_{s+1},\cdots,{\mathcal{C}}_{i-1},\tilde{\mathcal{C}}_1, \cdots,\tilde{\mathcal{C}}_{k+1},{\mathcal{C}}_{i+1},\cdots,{\mathcal{C}}_{n+1}\rangle }$\\
  =  $\mathcal{C}\bigcap  { \langle  {\mathcal{C}}_1,\cdots,{\mathcal{C}}_{s-1} \rangle }^{\bot}\bigcap  \leftidx{^{\bot}}  { \langle  {\mathcal{C}}_{s+1},\cdots,{\mathcal{C}}_{i-1},{\mathcal{C}}_{i},{\mathcal{C}}_{i+1},\cdots,{\mathcal{C}}_{n+1}\rangle }$\\  =  ${\mathcal{C}}_s$.

 similarly, if  ${\mathcal{C}}_s \in {\{{\mathcal{C}}_{i+1},\cdots,{\mathcal{C}}_{n+1}\}}$,then

 $\mathcal{C}\bigcap    \langle  {\mathcal{C}}_{1},\cdots,{\mathcal{C}}_{i-1},\tilde{\mathcal{C}}_1, \cdots,\tilde{\mathcal{C}}_{k+1},{\mathcal{C}}_{i+1},\cdots,{\mathcal{C}}_{s-1} \rangle^{\bot} \bigcap  \leftidx{^{\bot}}{ \langle  {\mathcal{C}}_{s+1},\cdots,{\mathcal{C}}_{n+1} \rangle }$\\
  = ${\mathcal{C}}_s$.

If $\tilde {\mathcal{C}}_s \in {\{\tilde{\mathcal{C}}_1, \tilde{\mathcal{C}}_2,\cdots,\tilde{\mathcal{C}}_{k+1}\}}$, then

 $\mathcal{C}\bigcap    \langle  {\mathcal{C}}_{1},\cdots,{\mathcal{C}}_{i-1},\tilde{\mathcal{C}}_1, \cdots,\tilde{\mathcal{C}}_{s-1} \rangle ^{\bot} \bigcap  \leftidx{^{\bot}} { \langle  \tilde{\mathcal{C}}_{s+1},\cdots,\tilde{\mathcal{C}}_{k+1},\cdots,{\mathcal{C}}_{i+1},\cdots,{\mathcal{C}}_{n+1} \rangle }$\\
  =  $\mathcal{C}\bigcap    \langle  {\mathcal{C}}_{1},\cdots,{\mathcal{C}}_{i-1}\rangle ^{\bot}\bigcap \langle \tilde{\mathcal{C}}_1, \cdots,\tilde{\mathcal{C}}_{s-1} \rangle ^{\bot} \bigcap  \leftidx{^{\bot}} { \langle  \tilde{\mathcal{C}}_{s+1},\cdots,\tilde{\mathcal{C}}_{k+1}\rangle } \bigcap  \leftidx{^{\bot}} { \langle  {\mathcal{C}}_{i+1},\cdots,}
 {\mathcal{C}}_{n+1}\rangle $\\
  = ${\mathcal{C}}_i\bigcap \langle \tilde{\mathcal{C}}_1, \cdots,\tilde{\mathcal{C}}_{s-1} \rangle ^{\bot} \bigcap  \leftidx{^{\bot}} { \langle  \tilde{\mathcal{C}}_{s+1},\cdots,\tilde{\mathcal{C}}_{k+1}\rangle }$\\
  = $\tilde{\mathcal{C}}_s$.

Thus, the checking of the first condition of definition 2.3 is finished.

Step 2. Without losing of generality, we may assume $1\leq s\leq k$, and we want to check
$(\langle  {\mathcal{C}}_1,\cdots,\tilde {\mathcal{C}}_{s} \rangle ,\langle  \tilde {\mathcal{C}}_{s+1},\cdots,{\mathcal{C}}_{n+1} \rangle )$ is a \text{torsion pair} on $\mathcal{C}$.

Given  $X\in \mathcal{C}$, because $(\langle  {\mathcal{C}}_1,\cdots,{\mathcal{C}}_{i-1} \rangle ,\langle  {\mathcal{C}}_{i},\cdots,{\mathcal{C}}_{n+1} \rangle )$ is a \text{torsion pair} on $\mathcal{C}$, there is an exact sequence
\[ \begin{CD}
0 @>>> X_1 @> i_1>> X @> \pi_1>>  X_2 @>>>0
\end{CD} \]
such that $X_1 \in \langle  {\mathcal{C}}_1,\cdots,{\mathcal{C}}_{i-1} \rangle $ and $X_{2}\in \langle  {\mathcal{C}}_{i},\cdots,{\mathcal{C}}_{n+1} \rangle $.

By $\text{torsion pair}$ $(\langle  {\mathcal{C}}_1,\cdots,{\mathcal{C}}_{i} \rangle ,\langle  {\mathcal{C}}_{i+1},\cdots,{\mathcal{C}}_{n+1} \rangle )$,
there is an exact sequence
\[ \begin{CD}
0 @>>> X_3 @> i_2>>  X_2 @> \pi_2>>  X_4 @>>>0
\end{CD} \]
such that $X_3 \in \langle  {\mathcal{C}}_1,\cdots,{\mathcal{C}}_{i} \rangle $ and $X_{4}\in \langle  {\mathcal{C}}_{i+1},\cdots,{\mathcal{C}}_{n+1} \rangle $.

Because $X_3 \in { \langle  {\mathcal{C}}_1,\cdots,{\mathcal{C}}_{i-1} \rangle }^{\bot}$ since  $X_2 \in { \langle  {\mathcal{C}}_1,\cdots,{\mathcal{C}}_{i-1} \rangle }^{\bot}$, so $X_{3}\in  \mathcal{C}\bigcap \linebreak { \langle  {\mathcal{C}}_1,\cdots,{\mathcal{C}}_{i-1} \rangle }^{\bot}\bigcap  \leftidx{^{\bot}}  { \langle  {\mathcal{C}}_{i+1},\cdots,{\mathcal{C}}_{n+1} \rangle } = {\mathcal{C}}_i $.

By $\text{torsion pair} (\langle \tilde{\mathcal{C}}_1, \cdots,\tilde{\mathcal{C}}_{s}\rangle ,\langle \tilde{\mathcal{C}}_{s+1}, \cdots,\tilde{\mathcal{C}}_{k+1}\rangle )$ on $\mathcal{C}_i$,
there is an exact sequence
\[ \begin{CD}
0 @>>> X_5 @> i_3>>  X_3 @> \pi_3>>  X_6 @>>>0
\end{CD} \]
such that $X_5 \in \langle \tilde{\mathcal{C}}_1, \cdots,\tilde{\mathcal{C}}_{s}\rangle $ and $X_{6}\in \langle \tilde{\mathcal{C}}_{s+1}, \cdots,\tilde{\mathcal{C}}_{k+1}\rangle $.

By pushout of $i_2$ and $\pi_3$, we have the following commutative diagram£»
\[\begin{CD}
@. 0 @. 0\\
@. @VVV @VVV\\
@.X_5 @>>> X^{'}_5\\
@. @Vi_{3}VV @VVi_{4}V \\
0 @>>> X_3 @> i_2>>  X_2 @> \pi_2>>  X_4 @>>> 0\\
@. @V\pi_3VV @VV\pi_4V @|\\
0 @>>> X_6 @>>>X_\mathcal{F} @>>>X_4 @>>>0\\
@. @VVV @VVV \\
@.0 @. 0
\end{CD}\]

By snake lemma, $X_5 = X^{'}_5$, so we have an exact sequence \[\begin{CD}0 @>>> X_5 @> i_4>> X_2 @> \pi_4>>  X_\mathcal{F} @>>> 0\end{CD}\]
such that $X_\mathcal{F} \in \langle \tilde{\mathcal{C}}_{s+1}, \cdots,\tilde{\mathcal{C}}_{k+1},{\mathcal{C}}_{i+1},\cdots,{\mathcal{C}}_{n+1}\rangle$.

By pullback of $i_4$ and $\pi_1$, we have the following commutative diagram:

\[\begin{CD}
@.  @. 0 @. 0\\
@. @. @VVV @VVV \\
0 @>>>X_1 @>>>X_{\mathcal{T}} @>>>X_5@>>>0\\
@. @| @VVV @VVi_4V \\
0 @>>>X_{1} @>> i_1> X @>> \pi_1> X_{2}@>>>0\\
@. @. @VVV @VV\pi_4V\\
@. @. X^{'}_{\mathcal{F}}@>>> X_{\mathcal{F}}\\
@. @. @VVV @VVV\\
@. @. 0@. 0
\end{CD}\]

By snake lemma, $ X^{'}_{\mathcal{F}} =  X_{\mathcal{F}}$, so we have the following exact sequence:
 \[\begin{CD}0 @>>> X_\mathcal{T} @>>> X @>>> X_\mathcal{F} @>>> 0\end{CD}\]
 such that $ X_\mathcal{T}\in \langle {\mathcal{C}}_{1}, \cdots,\tilde{\mathcal{C}}_{s}\rangle $ and $X_\mathcal{F} \in \langle \tilde{\mathcal{C}}_{s+1}, \cdots,{\mathcal{C}}_{n+1}\rangle $.\\

Now we give the following definition which is very important to learn the structure of $n-\text{torsion pair}$.

\begin{defn} Series $\{({{\mathcal{T}}
_1,{\mathcal{F}}_1}),({{\mathcal{T}}
_2,{\mathcal{F}}_2}),\cdots,({{\mathcal{T}}
_n,{\mathcal{F}}_n})\}$ of \text{torsion pair}s on $\mathcal{C}$ is called an $n-$\text{torsion pair} series if
${\mathcal{T}}_1\subseteq{\mathcal{T}}_2\subseteq \cdots\subseteq{\mathcal{T}}_n$$
(\text{equivalently, }{\mathcal{F}}_1\supseteq {\mathcal{F}}_2 \supseteq \cdots \supseteq {\mathcal{F}}_n )$.
\end{defn}

The following  definition is an operation.

 \begin{defn}Let $(\mathcal{T},\mathcal{F})$ be a \text{torsion pair} on  $\mathcal{C}$, and  $\mathcal{D}$ be a subcategory of  $\mathcal{C}$. We call $(D^1_{(\mathcal{T},\mathcal{F})}(\mathcal{D}),D^2_{(\mathcal{T},\mathcal{F})}(\mathcal{D}))$ is a decomposition of $\mathcal{D}$ along $(\mathcal{T},\mathcal{F})$, where $D^1_{(\mathcal{T},\mathcal{F})}(\mathcal{D}) = \{X \mid $ There exists an exact sequence
$0\rightarrow X\rightarrow M\rightarrow Y\rightarrow 0$ such that  $X\in\mathcal{T}, Y\in\mathcal{F}, M\in\mathcal{D}\}$, $D^2_{(\mathcal{T},\mathcal{F})}(\mathcal{D}) = \{Y\mid$ There exists an exact sequence
$0\rightarrow X\rightarrow M\rightarrow Y\rightarrow 0$ such that  $X\in\mathcal{T}, Y\in\mathcal{F},  M\in\mathcal{D}$\}.
\end{defn}

 \begin{lem}If $\{({{\mathcal{T}}_1,{\mathcal{F}}_1}),({{\mathcal{T}}_2,{\mathcal{F}}_2})\}$ is a $2-\text{torsion pair}$ series on $\mathcal{C}$. Then $$\mathcal{F}_1 \bigcap \mathcal{T}_2 = D^2_{(\mathcal{T}_1,\mathcal{F}_1)}({\mathcal{T}}_2) = D^1_{(\mathcal{T}_2,\mathcal{F}_2)}({\mathcal{F}}_1).$$
 \end{lem}

\noindent{\bf Proof:} $\mathcal{F}_1 \bigcap \mathcal{T}_2 \subseteq D^2_{(\mathcal{T}_1,\mathcal{F}_1)}({\mathcal{T}}_2)$ is clear.

Suppose $X\in {\mathcal{T}}_2 $, by \text{torsion pair} $({{\mathcal{T}}_1,{\mathcal{F}}_1})$,there is an exact sequence
 $$ 0\longrightarrow X_{\mathcal{T}_1} \longrightarrow X \longrightarrow X_{\mathcal{F}_1}\longrightarrow 0$$
such that $X_{\mathcal{T}_1}\in \mathcal{T}_1 $ and $X_{\mathcal{F}_1}\in {\mathcal{F}_1}.$

However, $X_{\mathcal{F}_1}\in \leftidx{^\bot} \!\mathcal{F}_2$ since $X\in \leftidx{^\bot} \!\mathcal{F}_2$.
Thus,  $X_{\mathcal{F}_1}\in \leftidx{^\bot} \!\mathcal{F}_2 \bigcap \mathcal{C} =  {\mathcal{T}_2}$ and  $X_{\mathcal{F}_1}\in {\mathcal{T}_2}\bigcap {\mathcal{F}_1}$. So $\mathcal{F}_1 \bigcap \mathcal{T}_2 = D^2_{(\mathcal{T}_1,\mathcal{F}_1)}({\mathcal{T}}_2)$.

The other half is similar.
\\

$n-\text{torsion pair}$ series will give a filtration for every module which is demonstrated below.

\begin{prop} If $\{({{\mathcal{T}}
_1,{\mathcal{F}}_1}),({{\mathcal{T}}
_2,{\mathcal{F}}_2}),\cdots,({{\mathcal{T}}
_n,{\mathcal{F}}_n})\}$ is an $n-\text{torsion pair}$ seires on $\mathcal{C}$. Then for every module $X$ in  $\mathcal{C}$, there is a filtration:

\[\xymatrix{0 \ar@{ = }[0,1]&X_0 \ar[0,1]  &X_1\ar[0,1] \ar[1,-1] &\cdots\ar[0,1]&X_{n + 1}\ar[1,-1]\ar@{ = }[0,1]&X\\
            &S_1   && S_{n +1}}\]
such that $0 \ra X_i \ra X_{i + 1} \ra S_{i
+ 1} \ra 0$ is an exact sequence for $i  =  1, 2, \dots, n + 1$, and $S_1 \in \mathcal{T}_1, S_i \in \mathcal{F}_{i-1} \bigcap \mathcal{T}_i$ for $ 1 <i <  n + 1 $,
 $S_{n+1}\in \mathcal{F}_n$ and $X_j \in \mathcal{T}_j$ for $j <  n + 1$.
 \end{prop}

 \noindent {\bf Proof}. Using induction on $n$.

 $n = 1$, by the second condition of definition 2.1, there is a filtration
 \[\xymatrix{0 \ar@{ = }[0,1]&X_0 \ar[0,1]  &X_1\ar[0,1] \ar[1,-1] &X_{2}\ar[1,-1]\ar@{ = }[0,1]&X\\
            &S_1   & S_{2}}\]
such that $0 \ra X_1 \ra X_2 \ra S_2 \ra 0$ is an exact sequence and $X_1 \in \mathcal{T}_1,
 S_2\in \mathcal{F}_1$.

 Suppose that the proposition is true for $n = k$, let us consider $n = k+1$.
 By $\text{torsion pair}$ $({{\mathcal{T}}
_{k+1},{\mathcal{F}}_{k+1}})$ on $\mathcal{C}$,there is an exact sequence
$$ 0\longrightarrow X_{k+1} \longrightarrow X \longrightarrow S_{k+2}\longrightarrow 0$$ such that
$X_{k+1}\in {\mathcal{T}}_{k+1}$ and $S_{k+2}\in {\mathcal{F}}_{k+1}$.

Because  $\{({{\mathcal{T}}
_1,{\mathcal{F}}_1}),({{\mathcal{T}}
_2,{\mathcal{F}}_2}),\cdots,({{\mathcal{T}}
_k,{\mathcal{F}}_k})\}$ is a $k-\text{torsion pair}$ series on $\mathcal{C}$, by induction, there is a filtration:
\[\xymatrix{0 \ar@{ = }[0,1]&X_0 \ar[0,1] &X_1\ar[0,1] \ar[1,-1] &\cdots\ar[0,1]&X_{k}\ar[1,-1]\ar[0,1]&X_{k+1}\ar[1,-1] \ar@{ = }[0, 1] &X_{k+1}
\\
            &S_1  & & S_{k}&S_{k+1}}\]
 such that $0 \ra X_i \ra X_{i + 1} \ra S_{i
+ 1} \ra 0$ is an exact sequence for $i  =  1, 2, \dots, k + 1$, and $S_1 \in \mathcal{T}_1, S_i \in \mathcal{F}_{i-1} \bigcap \mathcal{T}_i $ for $1 < i <  k + 1 $,
$ S_{k+1}\in \mathcal{F}_k, X_i \in \mathcal{T}_i$ for all $i$.
 However $S_{k+1}\in \mathcal{T}_{k+1}$ since $X_{k+1}\in \mathcal{T}_{k+1}$. So $S_{k+1}\in \mathcal{F}_{k} \bigcap \mathcal{T}_{k+1}$. The filtration is given.\\

 \begin{prop} If $\{({{\mathcal{T}}
_1,{\mathcal{F}}_1}),({{\mathcal{T}}
_2,{\mathcal{F}}_2}),\cdots,({{\mathcal{T}}
_n,{\mathcal{F}}_n})\}$ is an $n-\text{torsion pair}$ series on $\mathcal{C}$.
Then $\mathcal{F}_{i} \bigcap \mathcal{T}_{i+k} = \langle \mathcal{F}_{i} \bigcap \mathcal{T}_{i+1},\mathcal{F}_{i+1} \bigcap \mathcal{T}_{i+2},\cdots,\mathcal{F}_{i+k-1} \bigcap \mathcal{T}_{i+k}\rangle $.
\end{prop}

\noindent{\bf Proof}.$"\supseteq"$ is obviously.\\
$"\subseteq": $For $X\in \mathcal{F}_{i} \bigcap \mathcal{T}_{i+k}$, by the above lemma, there is a filtration of $X$:

\[\xymatrix{0\ar@{ = }[0,1]&X_0\cdots\ar[0,1]&X_{i+1}\ar[1,-1]\ar[0,1]&
\cdots\ar[0,1] &X_{i+k}\ar[0,1]\ar[1,-1]&\cdots \ar[0,1]&X_{n+1}\ar[1,-1]\ar@{ = }[0,1]&X\\
               &S_{i+1}&&S_{i+k}&&S_{n+1}}\]
 such that $0 \ra X_i \ra X_{i + 1} \ra S_{i
+ 1} \ra 0$ is an exact sequence for $i  =  1, 2, \dots, n + 1$, and $S_1 \in \mathcal{T}_1, S_i \in \mathcal{F}_{i-1} \bigcap \mathcal{T}_i$ for $ 1 < i < n + 1$,
 $S_{n+1}\in \mathcal{F}_n, X_i \in \mathcal{T}_i$ for $i <  n + 1$.

 First, we claim that $X_0 = X_1 = \cdots = X_i = 0$.

 In fact, $\operatorname{Hom}(X_{i},X_{i+1}) = 0$ since $X_i\in \mathcal{T}_i$ and $X_{i+1}$ is submodule of $X$ belongs to $\mathcal{F}_i$. By the exact sequence
$0 \ra X_i \ra X_{i+1} \ra S_{i+1} \ra 0$, one gains $X_{i} = 0$. Hence $X_0 = X_1 = \cdots = X_{i-1} = 0$.

 Second, we claim that $X_{i+k+1} = X_{i+k+2} = \cdots = X_{n+1} = X$.

 In fact, $\operatorname{Hom}(X_{n+1},S_{n+1}) = 0$ since $X_{n+1} = X\in \mathcal{F}_{i} \bigcap \mathcal{T}_{i+k} $ and $S_{n+1}\in \mathcal{F}_{n}$. By exact sequence
$0 \ra X_n \ra  X_{n+1} \ra  S_{n+1} \ra  0$, one gains $S_{n+1} = 0$ and $X_n = X_{n+1} = X$. Similarly, we have $X_{i+k+1} = X_{i+k+2} = \cdots = X_{n-1} = X$.

 Now, we have the following filtration:

 \[\xymatrix{0 \ar@{ = }[0,1]&X_0 \ar[0,1]  &X_{i+1}\ar[0,1] \ar[1,-1]&\cdots\ar[0,1] &X_{i+k}\ar[1,-1]\ar[0,1]&X_{i+k+1}\ar[1,-1]\ar@{ = }[0,1]&X\\
            &S_{i+1}&   & S_{i+k}&S_{i+k+1}}\]
Thus $X\in \langle \mathcal{F}_{i} \bigcap \mathcal{T}_{i+1},\mathcal{F}_{i+1} \bigcap \mathcal{T}_{i+2},\cdots,\mathcal{F}_{i+k-1} \bigcap \mathcal{T}_{i+k}\rangle $.\\

 The following is the relation between $n-\text{torsion pair}$ and $n-\text{torsion pair} $ series.

 \begin{thm} There is a one to one correspondence between the set of $n-\text{torsion pair}$ series on $\mathcal{C}$ and the set of $n-\text{torsion pair}$ on $\mathcal{C}$:
 \[
 \begin{CD}
\left\{\begin{array}{c}({{\mathcal{T}}_1,{\mathcal{F}}_1}),\cdots,({{\mathcal{T}}_n,{\mathcal{F}}_n})\}
 :\\ n-\text{torsion pair series} \text{ on }\mathcal{C}\end{array} \right\} \autorightleftharpoons{$\alpha$}{$\beta$} \left \{\begin{array}{c}({\mathcal{C}}_1,{\mathcal{C}}_2,\cdots,{\mathcal{C}}_{n+1}): \\n-\text{torsion pair} \text{ on }\mathcal{C} \end{array} \right\}
 \end{CD}
 \]
such that
$\alpha(\{({{\mathcal{T}}
_1,{\mathcal{F}}_1}),({{\mathcal{T}}
_2,{\mathcal{F}}_2}),\cdots,({{\mathcal{T}}
_n,{\mathcal{F}}_n})\})  =  (\mathcal{T}
_1, {\mathcal{F}}_1\bigcap {\mathcal{T}
_2},\cdots,{\mathcal{F}}_{n-1}\bigcap {\mathcal{T}
_n},{\mathcal{F}
_n})$  and
$\beta(({\mathcal{C}}_1,{\mathcal{C}}_2,\cdots,{\mathcal{C}}_{n+1})) =
\{(\langle  {\mathcal{C}}_1,\cdots,{\mathcal{C}}_{i} \rangle ,\langle  {\mathcal{C}}_{i},\cdots,{\mathcal{C}}_{n+1} \rangle ) \mid {i = 1,2,\cdots,n}\}$.
\end{thm}

\noindent {\bf Proof:} First, we check that $(\mathcal{T}
_1, {\mathcal{F}}_1\bigcap {\mathcal{T}
_2},\cdots,{\mathcal{F}}_{n-1}\bigcap {\mathcal{T}
_n},{\mathcal{F}
_n})$ is an $n-\text{torsion pair}$ on $\mathcal{C}$ .

$\left( 1 \right){\mathcal{F}}_{i-1}\bigcap {\mathcal{T}_i} = \mathcal{C}\bigcap {\mathcal{T}_{i-1}}^\bot\bigcap\mathcal{C}\bigcap \leftidx{^\bot} \!\mathcal{F}_{i} =
\mathcal{C}\bigcap {\mathcal{T}_{i-1}}^\bot \bigcap \leftidx{^\bot} \!\mathcal{F}_{i} = \mathcal{C}\bigcap \langle \mathcal{T}
_1, {\mathcal{F}}_1\bigcap {\mathcal{T}
_2},\cdots, \linebreak {\mathcal{F}}_{i-2}\bigcap {\mathcal{T}
_{i-1}}\rangle ^\bot \bigcap \leftidx{^\bot} \!\langle {\mathcal{F}}_{i}\bigcap {\mathcal{T}_{i+1},\cdots,\mathcal{F}_{n}}\rangle $ by the above proposition.

$\left( 2 \right)$ Obviously, $(\langle \mathcal{T}
_1, {\mathcal{F}}_1\bigcap {\mathcal{T}
_2},\cdots,{\mathcal{F}}_{i-1}\bigcap {\mathcal{T}
_{i}}\rangle ,\langle {\mathcal{F}}_{i}\bigcap {\mathcal{T}_{i+1},\cdots,\mathcal{F}_{n}}\rangle )  =  (\mathcal{T}_{i},\mathcal{F}_{i})$.

Second, we check that $\{(\langle  {\mathcal{C}}_1,\cdots,{\mathcal{C}}_{i} \rangle ,\langle  {\mathcal{C}}_{i},\cdots,{\mathcal{C}}_{n+1} \rangle )\}_{i = 1,2,\cdots,n}$ is an $n-\text{torsion pair}$ series on $\mathcal{C}$. But this is clear.

Third, we check that $\beta\alpha = 1$.

$\beta\alpha(\{({{\mathcal{T}}
_1,{\mathcal{F}}_1}),({{\mathcal{T}}
_2,{\mathcal{F}}_2}),\cdots,({{\mathcal{T}}
_n,{\mathcal{F}}_n})\})  =  \beta(\mathcal{T}
_1, {\mathcal{F}}_1\bigcap {\mathcal{T}
_2},\cdots,{\mathcal{F}}_{n-1}\bigcap {\mathcal{T}
_n},{\mathcal{F}
_n})\\
 = \{({{\mathcal{T}}
_1,{\mathcal{F}}_1}),({{\mathcal{T}}
_2,{\mathcal{F}}_2}),\cdots,({{\mathcal{T}}
_n,{\mathcal{F}}_n})\}$  by the above proposition.

Last, we check that $\alpha\beta = 1$.

$\alpha\beta(({\mathcal{C}}_1,{\mathcal{C}}_2,\cdots,{\mathcal{C}}_{n+1})) = \alpha(
\{(\langle  {\mathcal{C}}_1,\cdots,{\mathcal{C}}_{i} \rangle ,\langle  {\mathcal{C}}_{i + 1},\cdots,{\mathcal{C}}_{n+1} \rangle )\}\mid {i = 1,2,\cdots,n})\\
 = \{
\langle  {\mathcal{C}}_{i},\cdots,{\mathcal{C}}_{n+1} \rangle \bigcap \langle  {\mathcal{C}}_1,\cdots,{\mathcal{C}}_{i} \rangle  \mid {i = 1,2,\cdots,n+1}\}\\
 = \{ \mathcal{C}\bigcap{ \langle  {\mathcal{C}}_1,\cdots,{\mathcal{C}}_{i-1} \rangle }^{\bot}\bigcap  \leftidx{^{\bot}}  { \langle  {\mathcal{C}}_{i+1},\cdots,{\mathcal{C}}_{n+1} \rangle } \mid {i = 1,2,\cdots,n+1}\}\\
 = ({\mathcal{C}}_1,{\mathcal{C}}_2,\cdots,{\mathcal{C}}_{n+1})$.

\begin{prop}$({\mathcal{C}}_1,{\mathcal{C}}_2,\cdots,{\mathcal{C}}_{n+1})$ is an $n-\text{torsion pair}$ on $\mathcal{C}$ if and only if\\
$\left( 1 \right)$ $\operatorname{Hom}(X,Y) = 0$ for all $X\in\mathcal{C}_i$, $Y\in\mathcal{C}_j, i < j$.\\
$\left( 2 \right)$ For every $X\in \mathcal{C}$, there is a filtration:

\[\xymatrix{0 \ar@{ = }[0,1]&X_0 \ar[0,1]  &X_1\ar[0,1] \ar[1,-1] &\cdots\ar[0,1]&X_{n + 1}\ar[1,-1]\ar@{ = }[0,1]&X\\
            &S_1   && S_{n +1}}\]
such that $0 \ra X_i \ra X_{i + 1} \ra S_{i
+ 1} \ra 0$  is an exact sequence and $S_i \in \mathcal{C}_i$ for all $i$.
\end{prop}

\noindent{\bf Proof:}
$"\Longrightarrow"$: Let $\mathcal{T}_{i} = \langle  {\mathcal{C}}_1,\cdots,{\mathcal{C}}_{i} \rangle ,\mathcal{F}_{i} = \langle  {\mathcal{C}}_{i},\cdots,{\mathcal{C}}_{n+1} \rangle , i = 1,2,\cdots,n.$ Then $\{({{\mathcal{T}}
_1,{\mathcal{F}}_1}),({{\mathcal{T}}
_2,{\mathcal{F}}_2}),\cdots,({{\mathcal{T}}
_n,{\mathcal{F}}_n})\}$ is an $n-\text{torsion pair}$ series by proposition 2.12.

By the proof of the above proposition, we know $\mathcal{C}_{i} = \mathcal{F}_{i-1}\bigcap \mathcal{T}_i$.

Hence, for every module $X$ in  $\mathcal{C}$, there is a filtration:

\[\xymatrix{0 \ar@{ = }[0,1]&X_0 \ar[0,1]  &X_1\ar[0,1] \ar[1,-1] &\cdots\ar[0,1]&X_{n + 1}\ar[1,-1]\ar@{ = }[0,1]&X\\
            &S_1   && S_{n +1}}\]
such that $0 \ra X_i \ra X_{i + 1} \ra S_{i
+ 1} \ra 0$ is an exact sequence  and $ S_i \in \mathcal{F}_{i-1} \bigcap \mathcal{T}_i = \mathcal{C}_i$.

$"\Longleftarrow"$: First, we show that ${\mathcal{C}}_i  =  \mathcal{C}\bigcap  { \langle  {\mathcal{C}}_1,\cdots,{\mathcal{C}}_{i-1} \rangle }^{\bot}\bigcap  \leftidx{^{\bot}}  { \langle  {\mathcal{C}}_{i+1},\cdots,{\mathcal{C}}_{n+1} \rangle }$ for $i = 1,2,\cdots,n+1$.

$"\subseteq"$ is clear;

$"\supseteq"$: $\forall X\in \mathcal{C}\bigcap  { \langle  {\mathcal{C}}_1,\cdots,{\mathcal{C}}_{i-1} \rangle }^{\bot}\bigcap  \leftidx{^{\bot}}  { \langle  {\mathcal{C}}_{i+1},\cdots,{\mathcal{C}}_{n+1} \rangle }$, there is a filtration:
\[\xymatrix{0 \ar@{ = }[0,1]&X_0 \ar[0,1] &X_1\ar[0,1] \ar[1,-1] &\cdots\ar[0,1]&X_{i-1}\ar[1,-1]\ar[0,1]&X_{i}\ar[1,-1]\ar[0,1]&
X_{i+1}\ar[1,-1]\ar[0,1]&\ar[0,1]\cdots \ar[0,1]&X_{n+1}\ar[1,-1]\ar@{ = }[0,1]&X\\
            &S_1  & & S_{i-1}&S_{i}&S_{i+1}&&S_{n+1}}\]

such that $0 \ra X_i \ra X_{i + 1} \ra S_{i
+ 1} \ra 0$ is an exact sequence  and $ S_i \in \mathcal{C}_i$.

Just like the proof of proposition 2.10, we have
$X_0 = X_1 = \cdots = X_{i-1} = 0$ and $X_i = X_{i+1} = \cdots = X_{n+1} = X$. So $X_{i} = S_{i}\in \mathcal{C}_i$.

Second, we show that
$(\langle  {\mathcal{C}}_1,\cdots,{\mathcal{C}}_{i} \rangle ,\langle  {\mathcal{C}}_{i+1},\cdots,{\mathcal{C}}_{n+1} \rangle )$ is a $\text{torsion pair}$ on $\mathcal{C}$ by definition 2.1:

(1) Clear!

(2)  $\forall X\in \mathcal{C}$, there is a filtration:
\[\xymatrix{0 \ar@{ = }[0,1]&X_0 \ar[0,1] &X_1\ar[0,1] \ar[1,-1] &\cdots\ar[0,1]&X_{i-1}\ar[1,-1]\ar[0,1]&X_{i}\ar[1,-1]\ar[0,1]&
X_{i+1}\ar[1,-1]\ar[0,1]&\ar[0,1]\cdots \ar[0,1]&X_{n+1}\ar[1,-1]\ar@{ = }[0,1]&X\\
            &S_1  & & S_{i-1}&S_{i}&S_{i+1}&&S_{n+1}}\]
such that $0 \ra X_i \ra X_{i + 1} \ra S_{i
+ 1} \ra 0$ is an exact sequence  and $ S_i \in \mathcal{C}_i$.

It is clear that $X_{i}\in \langle  {\mathcal{C}}_1,\cdots,{\mathcal{C}}_{i} \rangle $, we claim that $X/X_i \in \langle  {\mathcal{C}}_{i+1},\cdots,{\mathcal{C}}_{n+1} \rangle $.

In fact, by snake lemma we have the following commutative diagram:
\[\begin{CD}
@.  @. 0 @. 0\\
@. @. @VVV @VVV \\
0 @>>>X_i @>>>X_{i+1} @>>>S_{i+1}@>>>0\\
@. @| @VVV @VVV \\
0 @>>>X_i @>>>X_{i+2} @>>>X_{i+2}/X_{i}@>>>0\\
@. @. @VVV @VVV\\
@. @. S_{i+2} @ =  S_{i+2}\\
@. @. @VVV @VVV\\
@. @. 0@. 0
\end{CD}\]
Hence $X_{i+2}/X_{i}\in \langle  {\mathcal{C}}_{i+1},{\mathcal{C}}_{i+2} \rangle $.

Use snake lemma again,we have the following commutative diagram:
\[\begin{CD}
@.  @. 0 @. 0\\
@. @. @VVV @VVV \\
0 @>>>X_i @>>>X_{i+2} @>>>X_{i+2}/X_{i}@>>>0\\
@. @| @VVV @VVV \\
0 @>>>X_i @>>>X_{i+3} @>>>X_{i+3}/X_{i}@>>>0\\
@. @. @VVV @VVV\\
@. @. S_{i+3} @ =  S_{i+3}\\
@. @. @VVV @VVV\\
@. @. 0@. 0
\end{CD}\]
Hence $X_{i+3}/X_{i}\in \langle  {\mathcal{C}}_{i+1},{\mathcal{C}}_{i+2},{\mathcal{C}}_{i+3} \rangle $.

Similarly, we can obtain $X_{n+1}/X_i \in \langle  {\mathcal{C}}_{i+1},\cdots,{\mathcal{C}}_{n+1} \rangle $.

Now, $0 \ra X_i \ra X_{n+1} \ra X_{n+1}/X_{i} \ra 0 $ is the desired exact sequence.\\

The following lemma is well known\cite{10}.
\begin{lem} If $\mathcal{B}$ is a subcategory of $\Lambda-mod$, then $^ \bot \!((^ \bot \!\mathcal{B})^\bot) = ^ \bot \!\mathcal{B}$ and $(^ \bot \!(\mathcal{B}^\bot))^\bot = \mathcal{B}^\bot$ and $(^ \bot \!\mathcal{B},(^ \bot \!\mathcal{B})^\bot)$ and
$(\mathcal{B}^\bot,^ \bot \!(\mathcal{B}^\bot))$ are both \text{torsion pair}s.
\end{lem}

The following means that the condition $(2)$ in Definition 2.3 will be superfluous in some conditions.

\begin{cor} Let  ${\mathcal{C}}_1,{\mathcal{C}}_2,\cdots,{\mathcal{C}}_n$ be full subcategories of  $\Lambda\text{-mod}$,
 if
${\mathcal{C}}_i  =   \langle  {\mathcal{C}}_1,\cdots,$ \linebreak[4] ${\mathcal{C}}_{i-1} \rangle ^{\bot}\bigcap  \leftidx{^{\bot}}  { \langle  {\mathcal{C}}_{i+1},\cdots,{\mathcal{C}}_{n+1} \rangle }$ for $i = 1,2,\cdots,n+1$. Then $({\mathcal{C}}_1,{\mathcal{C}}_2,\cdots,{\mathcal{C}}_n)$ is an $n-\text{torsion pair}$ on $\Lambda\text{-mod}$.
\end{cor}
\noindent{\bf Proof:} It is enough to show the second condition of the above proposition since the first condition is clear.

 By the above lemma, there is a fact: $(\leftidx{^\bot} \!\mathcal{C}_{n+1},\mathcal{C}_{n+1})$ is a \text{torsion pair} since ${\mathcal{C}}_{n+1} =  { \langle  {\mathcal{C}}_1,\cdots,{\mathcal{C}}_{n} \rangle }^{\bot}$.

Now, we use induction on $n$ to show.

If $n = 1$, clear.

Suppose that the proposition is true for $n = k\geq 1$, we consider the case of $n = k+1$.

Step 1, claim:$\langle \mathcal{C}_{k+1},\mathcal{C}_{k+2}\rangle   =  { \langle  {\mathcal{C}}_1,\cdots,{\mathcal{C}}_{k} \rangle }^{\bot}$.

In fact, $"\subseteq"$ is clear.

$"\supseteq"$: $\forall X\in { \langle  {\mathcal{C}}_1,\cdots,{\mathcal{C}}_{k} \rangle }^{\bot}$, by \text{torsion pair} $(\leftidx{^\bot} \!\mathcal{C}_{k+2},\mathcal{C}_{k+2})$, $\exists$ an exact sequence  $0 \ra X_{k+1} \ra X \ra T_{k
+ 2} \ra 0$ such that $X_{k+1}\in \leftidx{^\bot} \!\mathcal{C}_{k+2} $ and $T_{k+ 2}\in \mathcal{C}_{k+2}$.
$X_{k+1}\in { \langle  {\mathcal{C}}_1,\cdots,{\mathcal{C}}_{k} \rangle }^{\bot}$ since  $X\in { \langle  {\mathcal{C}}_1,\cdots,{\mathcal{C}}_{k} \rangle }^{\bot}$. Thus $X_{k+1}\in \mathcal{C}_{k+1}$ and $X\in \langle \mathcal{C}_{k+1},\mathcal{C}_{k+2}\rangle $.

Step 2.  By induction, $(\mathcal{C}_1,\cdots,{\mathcal{C}}_{k},\langle \mathcal{C}_{k+1},\mathcal{C}_{k+2}\rangle )$ is a $k-\text{torsion pair}$ on $\Lambda\text{-mod}$. So $\forall X\in \Lambda\text{-mod}$, there is a filtration:
\[\xymatrix{0 \ar@{ = }[0,1]&X_0 \ar[0,1] &X_1\ar[0,1] \ar[1,-1] &\cdots\ar[0,1]&X_{k-1}\ar[1,-1]\ar[0,1]&X_{k}\ar[1,-1]\ar[0,1]&
X_{k+1}\ar[1,-1]\ar@{ = }[0,1]&X\\
            &S_1  & & S_{k-1}&S_{k}&S}\]
            such that $S_i\in \mathcal{C}_i$ and $S\in \langle \mathcal{C}_{k+1},\mathcal{C}_{k+2}\rangle $.

By \text{torsion pair} $(\leftidx{^\bot} \!\mathcal{C}_{k+2},\mathcal{C}_{k+2})$, there is an exact sequence  $0 \ra S_{k+1}\ra S \ra S_{k+ 2} \ra 0$ such that $S_{k+1}\in \leftidx{^\bot} \!\mathcal{C}_{k+2}$ and $S_{k+2}\in \mathcal{C}_{k+2}$.

Because $S\in \langle \mathcal{C}_{k+1},\mathcal{C}_{k+2}\rangle $, then $S\in \langle \mathcal{C}_1,\cdots,{\mathcal{C}}_{k}\rangle ^\bot$, so $S_{k+1}\in \langle \mathcal{C}_1,\cdots,{\mathcal{C}}_{k}\rangle ^\bot$, hence $S_{k+1}\in \mathcal{C}_{k+1}$ since $S_{k+1}\in \leftidx{^\bot} \!\mathcal{C}_{k+2}$.

By pullback of $(X\rightarrow S,S_{k+1}\rightarrow S)$, we have the following commutative diagram:

\[\begin{CD}
@.  @. 0 @. 0\\
@. @. @VVV @VVV \\
0 @>>>X_k @>>>X_{k+1} @>>>S_{k+1}@>>>0\\
@. @| @VVV @VVV \\
0 @>>>X_k @>>>X @>>>S@>>>0\\
@. @. @VVV @VVV\\
@. @. S_{k+2} @ =  S_{k+2}\\
@. @. @VVV @VVV\\
@. @. 0@. 0
\end{CD}\]

Now, we find a filtration:

\[\xymatrix{0 \ar@{ = }[0,1]&X_0 \ar[0,1] &X_1\ar[0,1] \ar[1,-1] &\cdots\ar[0,1]&X_{k-1}\ar[1,-1]\ar[0,1]&X_{k}\ar[1,-1]\ar[0,1]&
X_{k+1}\ar[1,-1]\ar[0,1]&X_{k+2}\ar[1,-1]\ar@{ = }[0,1]&X\\
            &S_1  & & S_{k-1}&S_{k}&S_{k+1}&S_{k+2}}\]
           such that $0 \ra X_i \ra X_{i + 1} \ra S_{i
+ 1} \ra 0$ is an exact sequence  and $ S_i \in \mathcal{C}_i$.
\begin{prop} Let $({\mathcal{C}}_1,{\mathcal{C}}_2,\cdots,{\mathcal{C}}_n)$ is an $n-\text{torsion pair}$ on $\mathcal{C}$. Then\\
$\left( 1 \right)$ $\langle  {\mathcal{C}}_{i+1},\cdots,{\mathcal{C}}_{i+k} \rangle    =  \mathcal{C}\bigcap  { \langle  {\mathcal{C}}_1,\cdots,{\mathcal{C}}_{i} \rangle }^{\bot}\bigcap  \leftidx{^{\bot}}  { \langle  {\mathcal{C}}_{i+k+1},\cdots,{\mathcal{C}}_{n+1} \rangle }$\\
$\left( 2 \right)$ $({\mathcal{C}}_i,{\mathcal{C}}_{i+1},\cdots,{\mathcal{C}}_{i+k})$ is a $k-\text{torsion pair}$ on $\langle {\mathcal{C}}_i,{\mathcal{C}}_{i+1},\cdots,{\mathcal{C}}_{i+k}\rangle $\\
$\left( 3 \right)$ $({\mathcal{C}}_1,\cdots,{\mathcal{C}}_{i-1},\langle {\mathcal{C}}_i,{\mathcal{C}}_{i+1},\cdots,{\mathcal{C}}_{i+k}\rangle , {\mathcal{C}}_{i+k+1},\cdots,{\mathcal{C}}_{n+1})$ is an $(n-k)-\text{torsion pair}$.
\end{prop}

\noindent {\bf Proof}.(1) $"\subseteq"$: clear!

 $"\supseteq"$: $\forall X\in\mathcal{C}\bigcap  { \langle  {\mathcal{C}}_1,\cdots,{\mathcal{C}}_{i} \rangle }^{\bot}\bigcap  \leftidx{^{\bot}}  { \langle  {\mathcal{C}}_{i+k+1},\cdots,{\mathcal{C}}_{n+1} \rangle }$, there is a filtration:
 \[\xymatrix{0\ar@{ = }[0,1]&X_0\cdots\ar[0,1]&X_{i+1}\ar[1,-1]\ar[0,1]&
\cdots\ar[0,1] &X_{i+k}\ar[0,1]\ar[1,-1]&\cdots \ar[0,1]&X_{n+1}\ar[1,-1]\ar@{ = }[0,1]&X\\
               &S_{i+1}&&S_{i+k}&&S_{n+1}}\]
               such that $X_0 = X_1 = \cdots = X_{i-1} = 0$ and $X_{i+k+1} = X_{i+k+2} = \cdots = X_{n+1} = X$.

(2) Checking by Definition $2.3$, the first condition holds by (1), and the second condition holds by similar techniques in proof of proposition 2.10 and (1).

(3) Checking by Definition $2.3$, the first condition obviously holds, $\forall X\in \mathcal{C}$, there is a filtration:
 \[\xymatrix{0\ar@{ = }[0,1]&X_0\cdots\ar[0,1]&X_{i+1}\ar[1,-1]\ar[0,1]&
\cdots\ar[0,1] &X_{i+k}\ar[0,1]\ar[1,-1]&\cdots \ar[0,1]&X_{n+1}\ar[1,-1]\ar@{ = }[0,1]&X\\
               &S_{i+1}&&S_{i+k}&&S_{n+1}}\]
              use the similar techniques in the last part of proof of proposition 2.12, we have the following exact sequence:
    \[\begin{CD} 0 @>>>X_i @>>>X_{i+k} @>>>X_{i+k}/X_{i}@>>>0\end{CD}\]

    Let $\hat S = X_{i+k}/X_{i}$, then we have the desired filtration:
    \[\xymatrix{0\ar@{ = }[0,1]&X_0\cdots\ar[0,1]&X_{i+1}\ar[1,-1]\ar[0,1]&X_{i+k}\ar[0,1]\ar[1,-1]&\cdots \ar[0,1]&X_{n+1}\ar[1,-1]\ar@{ = }[0,1]&X\\
               &S_{i+1}&\hat S&&S_{n+1}}\]

\begin{cor} Suppose $\{({{\mathcal{T}}
_1,{\mathcal{F}}_1}),({{\mathcal{T}}
_2,{\mathcal{F}}_2}),\cdots,({{\mathcal{T}}
_n,{\mathcal{F}}_n})\}$ is an $n-\text{torsion pair}$ series on $\mathcal{C}$. Let $\mathcal{F}_{0} = \mathcal{T}_{n+1} = \mathcal{C}, \mathcal{F}_{n+1} = \mathcal{T}_{0} = 0 $. Then $\{(\mathcal{T}_{i+1}\bigcap \mathcal{F}_{i},\mathcal{F}_{i+1}\bigcap \mathcal{T}_{i+k+1}),\cdots, \linebreak (\mathcal{T}_{i+k}\bigcap \mathcal{F}_{i}, \mathcal{F}_{i+k}\bigcap \mathcal{T}_{i+k+1})\}$ is a $k-\text{torsion pair}$ seriies on $\mathcal{T}_{i+k+1}\bigcap \mathcal{F}_{i}$ for $i = 0,1,\cdots, \linebreak n-1$ and $k > 0$.
\end{cor}

\noindent{\bf Proof}. Let $\mathcal{C}_1 = \mathcal{T}_1,\mathcal{C}_l = \mathcal{F}_{l-1}\bigcap \mathcal{T}_l(2\leq l \leq n),\mathcal{C}_{n+1} = \mathcal{F}_{n+1}$.
 So  $({\mathcal{C}}_1,{\mathcal{C}}_2,\cdots,{\mathcal{C}}_{n+1})$ is an $n-\text{torsion pair}$ on $\mathcal{C}$, and $({\mathcal{C}}_i,{\mathcal{C}}_{i+1},\cdots,{\mathcal{C}}_{i+k})$ is a $k-\text{torsion pair}$ on $\langle {\mathcal{C}}_i,{\mathcal{C}}_{i+1},\cdots,\linebreak {\mathcal{C}}_{i+k}\rangle$ . Thus $\{(\langle {\mathcal{C}}_i,\cdots,{\mathcal{C}}_{i+l}\rangle ,\langle {\mathcal{C}}_{i+l+1},\cdots,{\mathcal{C}}_{i+k+l}\rangle )\mid l = 1,2,\cdots,k\}$ is a $k-\text{torsion pair}$ series on $\langle {\mathcal{C}}_i,{\mathcal{C}}_{i+1},\cdots,{\mathcal{C}}_{i+k}\rangle $. But $\langle {\mathcal{C}}_i,\cdots,{\mathcal{C}}_{i+l}\rangle = \mathcal{F}_{i}\bigcap \mathcal{T}_{i+l}, \ \langle {\mathcal{C}}_{i+l+1},\cdots,{\mathcal{C}}_{i+k+l}\rangle  = \mathcal{F}_{i+1}\bigcap \mathcal{T}_{i+k+l}$. The corollary is proved.

 \begin{cor} If  $({\mathcal{D}}_1,{\mathcal{D}}_2,\cdots,{\mathcal{D}}_{n+1})$ is a
defect $n-\text{torsion pair}$ on $\mathcal{C}$. Then there is an unique $n-\text{torsion pair}$  $({\mathcal{C}}_1,{\mathcal{C}}_2,\cdots,{\mathcal{C}}_{n+1})$ on $\mathcal{C}$ such that $\mathcal{D}_i\subseteq \mathcal{C}_i$.
\end{cor}

\noindent {\bf Proof}. Let $\mathcal{T}_{i} = \langle  {\mathcal{D}}_1,\cdots,{\mathcal{D}}_{i} \rangle , \mathcal{F}_{i} = \langle  {\mathcal{D}}_{i+1},\cdots,{\mathcal{D}}_{n} \rangle  $, Then $\{({{\mathcal{T}}
_1,{\mathcal{F}}_1}),({{\mathcal{T}}
_2,{\mathcal{F}}_2}),\cdots, \linebreak ({{\mathcal{T}}
_n, {\mathcal{F}}_n})\}$ is an $n-\text{torsion pair}$ series on $\mathcal{C}$.

Let $\mathcal{C}_{i} = \mathcal{F}_{i-1}\bigcap \mathcal{T}_i$, then $({\mathcal{C}}_1,{\mathcal{C}}_2,\cdots,{\mathcal{C}}_{n+1})$ is an $n-\text{torsion pair}$ on $\mathcal{C}$ such that $\mathcal{D}_i\subseteq \mathcal{C}_i$.

Suppose  $({\mathcal{C}}^{'}_1,{\mathcal{C}}^{'}_2,\cdots,{\mathcal{C}}^{'}_{n+1})$ is an other $n-\text{torsion pair}$ on $\mathcal{C}$ such that $\mathcal{D}_i\subseteq \mathcal{C}^{'}_i$, then $\mathcal{T}_{i} = \langle  {\mathcal{D}}_1,\cdots,{\mathcal{D}}_{i} \rangle \subseteq \langle  {\mathcal{C}}^{'}_1,\cdots,{\mathcal{C}}^{'}_{i} \rangle  = \mathcal{T}^{'}$, Similarly, $\mathcal{F}\subseteq\mathcal{F}^{'}$.
Therefore, $\mathcal{C}_{i} = \mathcal{F}_{i-1}\bigcap \mathcal{T}_i = \mathcal{F}^{'}_{i-1}\bigcap \mathcal{T}^{'}_i = \mathcal{C}^{'}_{i}$.\\

The following proposition is very useful.
\begin{prop} Suppose $\{({{\mathcal{T}}_1,   {\mathcal{F}}_1}),   ({{\mathcal{T}}_2,   {\mathcal{F}}_2})\}$ is a $2-\text{torsion pair}$ series on $\mathcal{C}$. Then
 we have the following 1 to 1 correspondence :\\
  \{$(\mathcal{T^\prime},   \mathcal{F^\prime})$: $1-\text{torsion pair}$ on ${\mathcal{F}}_1\bigcap {\mathcal{T}}_2$ \} $\autorightleftharpoons{F}{G}$ \{$({{\mathcal{T}}_3,  {\mathcal{F}}_3})$: $1-\text{torsion pair}$ on $\mathcal{C}$ such that ${\mathcal{T}}_1 \subseteq {\mathcal{T}}_3 \subseteq {\mathcal{T}}_2$\}
  where $F((\mathcal{T^\prime},   \mathcal{F^\prime}))  =  (\langle {{\mathcal{T}}_1}, \mathcal{T^\prime}\rangle ,   \langle \mathcal{F^\prime}, {\mathcal{F}}_2\rangle ), G(({\mathcal{T}}_3,  {\mathcal{F}}_3)) =  ({\mathcal{T}}_3 \bigcap {{\mathcal{F}}_1,  {\mathcal{F}}_3}\bigcap {\mathcal{T}}_2)$.
\end{prop}

\noindent{\bf Proof:} By Proposition 2.5,Theorem $2.11$ and Proposition 2.15, it is clear.

\begin{rem} The above lemma has a lot of generalized forms since we have so many results. And those forms can give a finer characterization for \text{torsion pair}s and module categories. For example, Theorem $2.1$ in  $\cite{1}$.
\end{rem}

The following is an example of $n-\text{torsion pair}$.
\begin{exm} Let $T$ be a tilting module, ${T_1, T_2,
\dots, T_n}$ be all non-isomorphic indecomposable summands of $T$. Then
$$(Gen(T_1), T_1^\bot \bigcap Gen(T_1 \oplus T_2), \dots, (T_1
\oplus \dots \oplus T_{n - 1})^\bot \bigcap Gen(T))$$
is an $(n-1)-\text{torsion pair}$ on  $ Gen(T)$.
\end{exm}

\noindent{\bf Proof}. Let$X_i = T_1\oplus \cdots \oplus T_i$, then $(GenX_i,X_i^\bot)$ is a \text{torsion pair} on $\Lambda\text{-mod}$. And
$\{(GenX_i,X_i^\bot)\mid i = ,1,2,\cdots,n\}$ is an $n-\text{torsion pair}$ series on $\Lambda\text{-mod}$. Therefore, $$(Gen(T_1), T_1^\bot \bigcap Gen(T_1 \oplus T_2), \dots, (T_1
\oplus \dots \oplus T_{n - 1})^\bot \bigcap Gen(T_n), T^\bot)$$
is an $n-\text{torsion pair}$ on $\Lambda\text{-mod}$ byTheorem $2.11$. So $(Gen(T_1), T_1^\bot \bigcap Gen(T_1 \oplus T_2), \dots, (T_1
\oplus \dots \oplus T_{n - 1})^\bot \bigcap Gen(T))$
is an $(n-1)-\text{torsion pair}$ on  $ Gen(T)$ by Proposition $2.15$.

\section{Decomposition by projective and injective modules}
 In this section, we always suppose $\A$ is an artin algebra. For given artin algebra $\Gamma$, we denote: $\mathcal{P}(\Gamma)$ is the category of all projective modules in $\Gamma\text{-mod}$, $\mathcal{I}(\Gamma)$ is the category of all injective modules in $\Gamma\text{-mod}$;  $\mathbf{E}(\Gamma) = \{(\mathcal{T},\mathcal{F})$ is \text{torsion pair} on $\Gamma\text{-mod}$ $\mid \mathcal{T}\bigcap \mathcal{P}(\Gamma) = \mathcal{F}\bigcap \mathcal{I}(\Gamma) = \phi\}$.  For a set $\Psi$ we denote the number of the elements of $\Psi$ by $\# \Psi$. For a subcategory $\mathcal{D}$ of $\A$-mod, let $\operatorname{Ind}\mathcal{D}$ be the set of pairwise non-isomorphic indecomposable modules in $\mathcal{D}$. For a module $M$, let $\operatorname{Ind}M$  =  $\operatorname{Ind} (\operatorname{add}M)$

\begin{defn} Suppose $\mathcal{C}$ is a full subcategory of $\A$-mod. A $\A$-module $M$ is called $\Ext$-projective in $\mathcal{C}$ if $\Ext^1_\Lambda(M,\mathcal{C}) = 0$. Dually, it is called $\Ext$-injective in $\mathcal{C}$ if $\Ext^1_\Lambda(\mathcal{C},M) = 0$.
\end{defn}

The following lemma is from \cite{1}.
\begin{lem}$\left( 1 \right)$ $(\Lambda e)^\bot  =  \leftidx{^\bot}(D(e\Lambda)) = \Lambda/\Lambda e\Lambda\text{-mod}$\\
$\left( 2 \right)$ $(Gen(\Lambda e),\Lambda/\Lambda e\Lambda\text{-mod})$ and $(\Lambda/\Lambda e\Lambda\text{-mod},Cogen(D(e\Lambda)))$ are both \text{torsion pair}s on $\Lambda\text{-mod}$.
\end{lem}

\noindent {\bf Proof:} It is clear that $\Lambda/\Lambda e\Lambda\text{-mod} = \{M\in \Lambda\text{-mod} \mid\ eM = 0\}$.

We claim: $(\Lambda e)^\bot = \{M\in \Lambda\text{-mod} \mid\ eM = 0\} = \leftidx{^\bot}(D(e\Lambda))$.

In fact, for any $M\in (\Lambda e)^\bot, \operatorname{Hom}(\Lambda e,M) = eM = 0$; For any $M\in \leftidx{^\bot}(D(e\Lambda)), \operatorname{Hom}\linebreak (M,D(e\Lambda)) = \operatorname{Hom}(M,\operatorname{Hom}(e\Lambda,J)) = \operatorname{Hom}(e\Lambda\otimes M,J) = D(eM) = 0$ $\Longleftrightarrow eM = 0$.

By (1), (2) is clear.

\begin{lem} Let $(\mathcal{C}_1,\mathcal{C}_2,\mathcal{C}_3)$ be a $2-\text{torsion pair}$ on $\Lambda\text{-mod}$:\\
$\left( 1 \right)$If $X\in \mathcal{C}_2$ is $\Ext$-projective in $\langle \mathcal{C}_2,\mathcal{C}_3\rangle $, $P_X\twoheadrightarrow X$ is the projective cover of $X$.
Then, there exists an exact sequence $0\rightarrow K_X\rightarrow P_X\rightarrow X\rightarrow 0$ such that $K_X\in \mathcal{C}_1$. Especially,
$P_X\in \langle \mathcal{C}_1,\mathcal{C}_2\rangle $ and $P_X\not\in\mathcal{C}_1$.\\
$\left( 2 \right)$ If $Y\in \mathcal{C}_2$ is $\Ext$-injective in $\langle \mathcal{C}_1,\mathcal{C}_2\rangle $, $Y\hookrightarrow I_Y$ is the injective envelope of $Y$.
Then, there exists an exact sequence $0\rightarrow Y\rightarrow I_Y\rightarrow C_Y\rightarrow 0$ such that $C_Y\in \mathcal{C}_3$. Especially,
$I_Y\in \langle \mathcal{C}_2,\mathcal{C}_3\rangle $ and $I_Y\not\in \mathcal{C}_3$.
\end{lem}

\noindent{\bf Proof:} We only proof (1); The proof of (2) is similar.

By $(\mathcal{C}_1, \langle \mathcal{C}_2,\mathcal{C}_3\rangle )$, there is a exact sequence $0\rightarrow K_X\rightarrow P_X\rightarrow L\rightarrow 0$ such that $K_X\in \mathcal{C}_1$ and $L\in \langle \mathcal{C}_2,\mathcal{C}_3\rangle $, obviously, there is an epimorphism $\eta:L\rightarrow X$ if we apply $\operatorname{Hom}(-, X)$ to the exact sequence.  Since $\operatorname{Ker}\eta \in \langle \mathcal{C}_2,\mathcal{C}_3\rangle $ and $X$ is $\Ext$-projective in $\langle \mathcal{C}_2,\mathcal{C}_3\rangle $,  $\eta$ is split. Thus $L = X\oplus \operatorname{Ker}\eta$, by the minimality of projective cover, $L = X$.

\begin{lem}  Let $(\mathcal{C}_1,\mathcal{C}_2,\mathcal{C}_3)$ be a $2-\text{torsion pair}$ on $\Lambda\text{-mod}$, and $X\in \Lambda\text{-mod}$ has a filtration
 \[\xymatrix{0 \ar@{ = }[0,1]&X_0 \ar[0,1]  &X_{1}\ar[0,1] \ar[1,-1] &X_{2}\ar[1,-1]\ar[0,1]&X_{3}\ar[1,-1]\ar@{ = }[0,1]&X\\
            &S_{1}   & S_{2}&S_{3}}\]
$\left( 1 \right)$ If $X$ is projective and $S_3 = 0$, then $S_2$ is $\Ext$-projective in $\langle \mathcal{C}_2,\mathcal{C}_3\rangle $ or $S_2 = 0$\\
$\left( 2 \right)$If $X$ is injective and $S_1 = 0$, then $S_2$ is $\Ext$-injective in $\langle \mathcal{C}_1,\mathcal{C}_2\rangle $ or $S_2 = 0$
\end{lem}

\noindent{\bf Proof:} We only proof (1); The proof of (2) is similar.

Since $S_3 = 0$, $X \cong X_3 \cong X_3$. Then $0 \ra X_1 \ra X \ra S_2 \ra 0$ is an exact sequence such that $X_1 \in \mathcal{C}_1, S_2 \in \langle \mathcal{C}_2,\mathcal{C}_3\rangle$. By Proposition 1.11 in Chapter 6 of \cite{13}, $S_2$ is $\Ext$-projective in $\langle \mathcal{C}_2,\mathcal{C}_3\rangle $.

\begin{prop} Let $({\mathcal{C}}_1,{\mathcal{C}}_2,\cdots,{\mathcal{C}}_{n+1})$ be an $n-\text{torsion pair}$ on $\A$-mod. Then there exists bijections:\\
$\left( 1 \right)$ $F: \operatorname{Ind} \ \mathcal{P}(\Lambda)\rightarrow \{X\in \operatorname{Ind} \ \mathcal{C}_i\mid$ $X$ is $\Ext$-projective in $ \langle {\mathcal{C}}_{i},{\mathcal{C}}_{i+1},\cdots,{\mathcal{C}}_{n+1}\rangle   \}$;\\
$\left( 2 \right)$ $G: \operatorname{Ind} \ \mathcal{I}(\Lambda)\rightarrow \{Y\in \operatorname{Ind} \ \mathcal{C}_j\mid$ $Y$ is $\Ext$-injective in $ \langle {\mathcal{C}}_{1},{\mathcal{C}}_{2},\cdots,{\mathcal{C}}_{j}\rangle   \}$.
\end{prop}

\noindent{\bf Proof:} We only proof (1); The proof of (2) is similar.\\
Step 1. For any indecomposable projective $\Lambda$-module $P$, there is a filtration
\[\xymatrix{0 \ar@{ = }[0,1]&X_0 \ar[0,1]  &\cdots\ar[0,1]&X_{i-1}\ar[1,-1]\ar[0,1]&X_{i}\ar[1,-1]\ar[0,1]&
X_{i+1}\ar[1,-1]\ar[0,1]&\ar[0,1]\cdots \ar[0,1]&X_{n+1}\ar[1,-1]\ar@{ = }[0,1]&P\\
              & & S_{i-1}&S_{i}&S_{i+1}&&S_{n+1}}\]
such that $0 \ra X_i \ra X_{i + 1} \ra S_{i
+ 1} \ra 0$ is an exact sequence  and $ S_i \in \mathcal{C}_i$.

Assume that $S_i\in\{S_1,S_2,\cdots,S_{n+1}\}$ is the last non-zero module, then $S_{i+1} = \cdots = S_{n+1} = 0$ and $X_i = X_{i+1} = \cdots = X_{n+1} = P$.

Now, we consider the following filtration
 \[\xymatrix{0 \ar@{ = }[0,1]&X_0 \ar[0,1]  &X^\prime_{i-1}\ar[0,1] \ar[1,-1] &X_{i}\ar[1,-1]\ar[0,1]&X_{i+1}\ar[1,-1]\ar@{ = }[0,1]&P\\
            &S^\prime_{i-1}   & S_{i}&S_{i+1}}\]

By lemma 3.4, $S_i$ is $\Ext$-projective in $ \langle {\mathcal{C}}_{i},{\mathcal{C}}_{i+1},\cdots,{\mathcal{C}}_{n+1}\rangle  $. We denote $F(P)  =  S_i$.\\
Step2. Suppose $X\in \operatorname{Ind} \mathcal{C}_i$ such that $X$ is $\Ext$-projective in $ \langle {\mathcal{C}}_{i},{\mathcal{C}}_{i+1},\cdots,{\mathcal{C}}_{n+1}\rangle$. Then we denote the projective cover of  $X$ By $P_X$ and denote $F^{-1}(X)  =  P_X$.\\
Step 3. It is clear that $F^{-1}F(P)  =  P$ for any indecomposable projective module $P$. On the other hand, since
$(\langle {\mathcal{C}}_{1}, \cdots,{\mathcal{C}}_{i - 1}\rangle , {\mathcal{C}}_{i},  \langle {\mathcal{C}}_{i+1}, \cdots,{\mathcal{C}}_{n+1}\rangle )$ is a $3-\text{torsion pair}$ on $\A$-mod, by Lemma 3.3, $FF^{-1}(X)  =  X$ for any $X\in \operatorname{Ind} \mathcal{C}_i$ which is $\Ext$-projective in $ \langle {\mathcal{C}}_{i},{\mathcal{C}}_{i+1},\cdots,{\mathcal{C}}_{n+1}\rangle $.

\begin{cor} Let $(\mathcal{T},\mathcal{F})$ be a \text{torsion pair} on $\Lambda\text{-mod}$. Then\\
$\left( 1 \right)$ there is an idempotent $e$ such that $\mathcal{T}\bigcap \mathcal{P}(\Lambda) = \operatorname{add}\Lambda e$, and $\mathcal{T}\bigcap (\Lambda e)^\bot$ has no $\Ext$-projective modules in $(\Lambda e)^\bot$;\\
$\left( 2 \right)$ there is an idempotent $e$ such that $\mathcal{F}\bigcap \mathcal{I}(\Lambda) = \operatorname{add} D(e\Lambda)$, and $\leftidx ^\bot  D(e\Lambda)\bigcap \mathcal{F}  $ has no $\Ext$-injective modules in $\leftidx ^\bot  D(e\Lambda)$.
\end{cor}

\noindent{\bf Proof:} We only proof (1); The proof of (2) is similar.

The first statement is clear, only the second one needs a proof:

$(Gen(\Lambda e),(\Lambda e)^\bot )$ is a \text{torsion pair} since $\Lambda e$ is a projective module. So we have a $2-\text{torsion pair}$ series $\{(Gen(\Lambda e),(\Lambda e)^\bot) ,(\mathcal{T},\mathcal{F})\}$, and we have a $2-\text{torsion pair}$ $(Gen(\Lambda e),(\Lambda e)^\bot \bigcap\mathcal{T},\mathcal{F})$.

Suppose that $X\in \mathcal{T}\bigcap (\Lambda e)^\bot$ is $\Ext$-projective in $(\Lambda e)^\bot$. Then obviously, $X\not\in Gen(\Lambda e)$. Let $f:P_X\twoheadrightarrow X$ is the projective cover of $X$. Then by proposition 3.4, $P_X \in \mathcal{T}$, and $X\in Gen(\Lambda e)$, this is a contradiction!

\begin{lem} Let $(\mathcal{C}_1,\mathcal{C}_2,\mathcal{C}_3)$ be a $2-\text{torsion pair}$ on $\Lambda\text{-mod}$:\\
$\left( 1 \right)$ If $\langle \mathcal{C}_1,\mathcal{C}_2\rangle $ is closed under kernel, $X\in \mathcal{C}_1$ is $\Ext$-projective in $\langle \mathcal{C}_1,\mathcal{C}_2\rangle $, and  $f:P_X\twoheadrightarrow X$ is the projective cover of $X$, then $P_X = X$ or $P_X\not\in\leftidx^\bot\mathcal{C}_3$;\\
$\left( 2 \right)$ If $\langle \mathcal{C}_2,\mathcal{C}_3\rangle $ is closed under cokernel, $X\in \mathcal{C}_3$ is $\Ext$-injective in $\langle \mathcal{C}_2,\mathcal{C}_3\rangle $, and  $g:X\hookrightarrow I_X$ is the injective envelope of $X$, then $I_X = X$ or $I_X\not\in{\mathcal{C}_1}^\bot$.
\end{lem}
\noindent{\bf Proof:} We only proof (1); The proof of (2) is similar.

Suppose $P_X\in\leftidx^\bot\mathcal{C}_3 = \langle \mathcal{C}_1,\mathcal{C}_2\rangle $, then exact sequence
\[ \begin{CD}
0 @>>> \operatorname{Ker} f @>>> P_X @> f>>  X @>>>0
\end{CD} \]
is split in $\langle \mathcal{C}_1,\mathcal{C}_2\rangle $ since $X\in \mathcal{C}_1$ is $\Ext$-projective in $\langle \mathcal{C}_1,\mathcal{C}_2\rangle $ and $Kerf \in \langle \mathcal{C}_1,\mathcal{C}_2\rangle $.

\begin{cor} Let $(\mathcal{T},\mathcal{F})$ be a \text{torsion pair} on $\Lambda\text{-mod}$.\\
$\left( 1 \right)$ If there are idempotents $e^0,e^1$ such that $\operatorname{add} \Lambda e^0\bigcap \operatorname{add} \Lambda e^1 = 0 $, $\mathcal{F}\bigcap \mathcal{I}(\Lambda)  =  \operatorname{add} D(e^0\Lambda)$, $\mathcal{T}\bigcap \mathcal{P}(\Lambda/\Lambda e^0\Lambda) = \operatorname{add} (\Lambda/\Lambda e^0\Lambda)e^1$. Then $\mathcal{T}\bigcap \mathcal{P}(\Lambda) = \phi$ if and only if for any $P\in \operatorname{add} \Lambda e^1$, $P\not \in \Lambda/\Lambda e^0\Lambda\text{-mod}$;\\
$\left( 2 \right)$ If there are orthogonal idempotents $\varepsilon^0,\varepsilon^1$ such that $\mathcal{T}\bigcap \mathcal{P}(\Lambda) = \operatorname{add} \Lambda e^0$, $\mathcal{F}\bigcap \mathcal{I}(\Lambda/\Lambda \varepsilon^0\Lambda) =  \operatorname{add} D(\varepsilon^1(\Lambda/\Lambda \varepsilon^0\Lambda))$. Then $\mathcal{F}\bigcap \mathcal{I}(\Lambda) = \phi$ if and only if for any $I \in \operatorname{add}\ D(\varepsilon^1\Lambda)$, $I\not \in \Lambda/\Lambda \varepsilon^0\Lambda\text{-mod}$.
\end{cor}

\noindent {\bf Proof:} We only proof (1); The proof of (2) is similar.\\
$"\Rightarrow"$  Since $(\A/\A e^0\A \text{-mod}, Cogen(D(e^0\A)))$ is a \text{torsion pair} by lemma 3.2, $(\mathcal{T}, \mathcal{F } \bigcap \A/\A e^0\A \text{-mod}, \linebreak Cogen(D(e^0\A))$ is a $2-\text{torsion pair}$ on $\A$-mod. Suppose $0 \neq P \in \operatorname{add} \Lambda e^1$. Then $P/e^0P \in \operatorname{add} (\Lambda/\Lambda e^0\Lambda)e^1$ and $P/e^0P \neq 0$. So by the above lemma,  $P  =  P/e^0P \in
\mathcal{T}$ or $P \not \in \A/\A e^0\A \text{-mod}$. Since $\mathcal{T}\bigcap \mathcal{P}(\Lambda) = \phi$, $P \not \in \A/\A e^0\A \text{-mod}$.\\
$"\Leftarrow"$ Suppose $\mathcal{T} \bigcap \mathcal{P}(\Lambda)\neq \phi$. Then there exists $0 \neq P \in \mathcal{T} \bigcap \mathcal{P}(\Lambda)$. Then $P$ is also projective in $\A/\A e^0\A \text{-mod}$. So $P \in \operatorname{add} \Lambda e^1$. This is a contradiction.\\

Now we start to show the structure of \text{torsion pair}s by decomposing them by projective modules and injective modules. First we give some notations .

We always assume that $ \Delta = \{e_1,e_2,\cdots,e_n\}$ is a fixed complete set of primitive orthogonal idempotents of $\Lambda$.
   Given $S = \{ \Delta_0,\Delta_1, \Delta_2,\cdots, \Delta_m\mid \Delta_i\subseteq \Delta \}$ such that $\Delta_1, \Delta_2,\cdots, \Delta_m\neq \phi$ and $\Delta_i\bigcap \Delta_j = \phi $ for $i\neq j$, we have the following notations :
$e_S^i = \sum_{e\in \Delta_i}e$, $\varepsilon_S^i = \sum^i_{j = 0} e_S^j$;
$\Lambda_S^0 = \Lambda,\Lambda_S^1 = \frac{\Lambda_S^0}{\Lambda_S^0e_S^0\Lambda_S^0} = \frac{\Lambda}{\Lambda\varepsilon_S^0\Lambda},\cdots,
\Lambda_S^{m + 1} =  \frac{\Lambda_S^{m}}{\Lambda_S^{m} e_S^m\Lambda_S^{m}}  = \frac{\Lambda}{\Lambda\varepsilon_S^m\Lambda}; \mathrm{P}_{i}(\Lambda_S^i) = \oplus_{e\in\Delta_i }\Lambda_S^i e, \mathrm{I}_{i}(\Lambda_S^i) = \oplus_{e\in\Delta_i }D(e\Lambda_S^i)$.

\begin{defn} Suoppose $S$ is  as the above. It is  called a 2-type part partition if:
$\left( 1 \right)$ $\forall 0< 2i\leq m$ and $e\in \Delta_{2i}$, $e_S^{2i-1}\Lambda_S^{2i-1}e\neq 0$;
$\left( 2 \right)$ $\forall 1< 2i+1\leq m$, and $e\in \Delta_{2i+1}$, $e\Lambda_S^{2i}e_S^{2i}\neq 0$.

Dually, $S$ is called a 2-type part partition if:
$\left( 1 \right)$ $\forall 0< 2i\leq m$ and $e\in \Delta_{2i}$, $e\Lambda_S^{2i-1}e_S^{2i-1}\neq 0$;
$\left( 2 \right)$ $\forall 1< 2i+1\leq m$, and $e\in \Delta_{2i+1}$, $e_S^{2i}\Lambda_S^{2i}e\neq 0$.
\end{defn}

\begin{lem} Let $I$ be an ideal of $\Lambda$, and $e,e^\prime$ be two idempotents. Then $\operatorname{Hom}_{\Lambda/I}((\Lambda/I)\cdot e, D(e^\prime \cdot\Lambda/I) = 0$ if and only if $e^\prime\cdot\Lambda/I\cdot e = 0$.
\end{lem}

\noindent {\bf Proof:} Notice that $\operatorname{Hom}_{\Lambda/I}((\Lambda/I)\cdot e, D(e^\prime \cdot\Lambda/I) = D(e^\prime \cdot\Lambda/I \cdot e)$.\\

We give the following notations for describing our theorem easily.

$\mathfrak{M} = \{(\mathcal{T},\mathcal{F})\mid  (\mathcal{T},\mathcal{F})$ is a \text{torsion pair} on $\Lambda\text{-mod}\}$;

$\mathfrak{N} = \{(S = \{ \Delta^\prime_0,\Delta^\prime_1, \Delta^\prime_2,\cdots, \Delta^\prime_m\},(\mathcal{T}^\prime,\mathcal{F}^\prime))\mid $ S is a 1-type part partition, $(\mathcal{T}^\prime,\mathcal{F}^\prime)\in \mathbf{E}(\Lambda_S^{m+1})\}$.

$\mathfrak{N}^\prime  =  \{(S = \{ \Delta^\prime_0,\Delta^\prime_1, \Delta^\prime_2,\cdots, \Delta^\prime_m\},(\mathcal{T}^\prime,\mathcal{F}^\prime)) \mid$ S is a 2-type part partition, $(\mathcal{T}^\prime,\mathcal{F}^\prime)\in \mathbf{E}(\Lambda_S^{m+1})\}$.\\

Now we are in a position to give a demonstration of how to decompose a \text{torsion pair} into $n-\text{torsion pair}$ by projective modules and injective modules.

Let$(\mathcal{T},\mathcal{F})$ be an \text{torsion pair} on $\Lambda\text{-mod}$:

a.Let $\mathcal{T}^0 = \mathcal{T},\mathcal{F}^0 = \mathcal{F},\Lambda^0 = \Lambda$, there exists some $\Delta_0\subseteq \Delta$ such that
$\mathcal{T}^0\bigcap \mathcal{P}(\Lambda^0) =  \operatorname{add} \oplus_{e\in\Delta_0 }\Lambda^0 e =   \operatorname{add}  \mathrm{P}_{0}(\Lambda^0)$. Let $\mathcal{T}^1 = \mathcal{T}^0
\bigcap (\mathrm{P_0}(\Lambda^0))^\bot, \mathcal{F}^1 = \mathcal{F}^0$, $\A^1  =  \A / \A e^0 \A $ where $e^0  =  \sum_{e\in \Delta_0}e$. Then $(\mathcal{T}^1, \mathcal{F}^1)$ is a \text{torsion pair} on $\Lambda^1\text{-mod}$ and $\mathcal{T}^1\bigcap \mathcal{P}(\Lambda^1) = \{ 0 \}$ by corollary 3.6. Hence we have a $2-\text{torsion pair}$ $(Gen\mathrm{P}_{0}(\Lambda^0),\mathcal{T}^1, \mathcal{F}^1)$ on $\Lambda\text{-mod}$;

b.There exists some $\Delta_1\subseteq \Delta - \Delta_0 $ such that
$\mathcal{F}^1\bigcap \mathcal{I}(\Lambda^1) =  \operatorname{add} \oplus_{e\in\Delta_1 }D(e\Lambda^1) =  \operatorname{add} \mathrm{I}_1(\Lambda^1)$. Let $\mathcal{T}^2 = \mathcal{T}^1, \mathcal{F}^2 = \mathcal{F}^1\bigcap \leftidx{^{\bot}} \! \mathrm{I}_1 (\Lambda^1)$, $\A^2  =  \A / \A \varepsilon^1 \A $ where $\varepsilon^1  =  \sum_{e\in \Delta_0 \bigcup \Delta_1}e$. Then $(\mathcal{T}^2, \mathcal{F}^2)$ is a \text{torsion pair} on $\Lambda^2$-mod and $\mathcal{F}^2\bigcap \mathcal{I}(\Lambda^2) = \{ 0 \}$ by corollary 3.6.  Hence we have a $3-\text{torsion pair}$ $(Gen \mathrm{P}_{0}(\Lambda^0),\mathcal{T}^2, \mathcal{F}^2,Cogen \mathrm{I}_1(\Lambda^1) )$ on $\Lambda$-mod;

The above operation  goes on alternatively, then it will eventually stop since $\#\Delta $ is finite.

Finally, we obtain:

(1) $\{ \Delta_0,\Delta_1, \Delta_2,\cdots, \Delta_m\mid \Delta_i\subseteq \Delta \}$ such that $\Delta_1, \Delta_2,\cdots, \Delta_m\neq \phi$ and $\Delta_i\bigcap \Delta_j = \phi $ for $i\neq j$;

(2) $(\mathcal{T}^{m+1}, \mathcal{F}^{m+1})$ is a \text{torsion pair} on $\Lambda^{m+1}-mod$ and $(\mathcal{T}^{m+1}, \mathcal{F}^{m+1})\in \mathbf{E}(\Lambda^{m+1})$;

(3) $(Gen \mathrm{P}_{0}(\Lambda^0),Gen \mathrm{P}_{2}(\Lambda^2),\cdots,\mathcal{T}^{m+1}, \mathcal{F}^{m+1},\cdots,Cogen \mathrm{I}_{3}(\Lambda^3),Cogen \mathrm{I}_{1}(\Lambda^1))$ is a $(m+2)-\text{torsion pair} $ on $\Lambda\text{-mod}$;

(4) $\Lambda = \Lambda^0\rightarrow \Lambda^1\rightarrow \cdots\rightarrow \Lambda^{m+1}$ is a series of quotient algebras.

\begin{thm} There is a  one to one  correspondence between $\mathfrak{M}$ and $\mathfrak{N}$:\\
$$\mathfrak{M} \autorightleftharpoons{F}{G} \mathfrak{N}$$

\end{thm}

\noindent {\bf Proof:} Step 1. Suppose $(\mathcal{T}, \mathcal{F}) \in \mathfrak{M}$. we use the above operation. Then we get $S = \{ \Delta_0,\Delta_1, \Delta_2,\cdots, \Delta_m\mid \Delta_i\subseteq \Delta \}$ and $(\mathcal{T}^{m+1}, \mathcal{F}^{m+1}) \in \mathbf{E}(\Lambda^{m+1})$, so we need to prove S is a 2-type part partition, but it follows from corollary 3.6 and lemma 3.8. Let $F((\mathcal{T}, \mathcal{F}))  =  (S, (\mathcal{T}^{m+1}, \mathcal{F}^{m+1}))$.\\
Step 2. Suppose $(S = \{ \Delta_0,\Delta_1, \Delta_2,\cdots, \Delta_m\},(\mathcal{T}^\prime,\mathcal{F}^\prime))\in \mathfrak{N}$. By induction on $m$. It is easy to see that $(Gen \mathrm{P}_{0}(\Lambda_S^0),Gen \mathrm{P}_{2}(\Lambda_S^2),\cdots,\mathcal{T}^\prime, \mathcal{F}^\prime,\cdots,Cogen \mathrm{I}_{3}(\Lambda_S^3),Cogen \mathrm{I}_{1}(\Lambda_S^1))$ is a $(m+2)-\text{torsion pair} $ on $\Lambda-mod$. Let $G((S,(\mathcal{T}^\prime, \mathcal{F}^\prime))) = (\mathcal{T}, \mathcal{F}) = (\langle Gen \mathrm{P}_{0}(\Lambda_S^0), \linebreak Gen \mathrm{P}_{2}(\Lambda_S^2),  \cdots,\mathcal{T}^\prime\rangle , \langle \mathcal{F}^\prime,\cdots, Cogen \mathrm{I}_{3}(\Lambda_S^3),Cogen \mathrm{I}_{1}(\Lambda_S^1)\rangle )$.

Claim:$\mathcal{T}\cap \mathcal{P}(\Lambda) = \operatorname{add}\ \mathrm{P}_0(\Lambda^0)$.

Otherwise, there exists some $e\in \Delta-\Delta_0$ such that $\Lambda e\in \mathcal{T}$. By proposition 3.5 and the above $(m+2)-\text{torsion pair}$, there exists $0\neq X\in Gen \mathrm{P}_{2i}(\Lambda_S^{2i})(\text{or}\ \mathcal{T}^\prime)$ for some $i\neq 0$, such that $X$ is $\Ext$-projective in $\Lambda_S^{2i}(\text{or}\ \Lambda_S^{m+1})$, and the projective cover  of $X$ is $\Lambda e$ since $\mathcal{T}^\prime\cap\mathcal{P}(\Lambda_S^{m+1}) = \phi $ and $X\in \mathrm{P}_{2i}(\Lambda_S^{2i})$. However, since $S$ is a 2-type part partition, $e_S^{2i-1}\Lambda_S^{2i-1}e\neq 0$. So Hom$_{\Lambda_S^{2i-1}}(X,D(e_S^{2i-1}\Lambda_S^{2i-1}))$ = Hom$_{\Lambda_S^{2i-1}}(\Lambda_S^{2i-1}e,D(e_S^{2i-1}\Lambda_S^{2i-1}))\neq 0$. Hence $X\not\in \leftidx^\bot\mathcal{F}$. So $\Lambda e\not\in \leftidx^\bot\mathcal{F}$. A contradiction!

Step by step, we know $F(\mathcal{T},\mathcal{F}) = (S,(\mathcal{T}^\prime,\mathcal{F}^\prime))$.

Step 3. Given $(\mathcal{T},\mathcal{F})\in \mathfrak{M}$, it is clear that $GF(\mathcal{T},\mathcal{F}) = (\mathcal{T},\mathcal{F})$.\\

Dually, if we start to decompose a \text{torsion pair} from the right hand (torsion-free class), Then we have the following theorem :
\begin{thm} There is a  one to one  correspondence between $\mathfrak{M}$ and $\mathfrak{N}^\prime$:\\
$$\mathfrak{M} \autorightleftharpoons{$F^\prime$}{G$^\prime$} \mathfrak{N}^\prime$$
\end{thm}

It's natural to ask that what is the relation between the above two kinds of decomposition.
The following theorem indicates that the decomposition of a \text{torsion pair} from left hand and right hand are the same.

\begin{thm} Suppose $(\mathcal{T},\mathcal{F})\in\mathfrak{M}$, $F((\mathcal{T},\mathcal{F})) = (S^\prime = \{ \Delta^\prime_0,\Delta^\prime_1, \Delta^\prime_2,\cdots, \Delta^\prime_u\},(\mathcal{T}^\prime,\mathcal{F}^\prime))$ and $F^\prime((\mathcal{T},\mathcal{F})) = (S^{\prime\prime} = \{ \Delta^{\prime\prime}_0,\Delta^{\prime\prime}_1, \Delta^{\prime\prime}_2,\cdots, \Delta^{\prime\prime}_v\},(\mathcal{T}^{\prime\prime},\mathcal{F}^{\prime\prime}))$. Then $(\mathcal{T}^{\prime},\mathcal{F}^{\prime}) = (\mathcal{T}^{\prime\prime},\mathcal{F}^{\prime\prime}).$
\end{thm}

\noindent {\bf Proof:} It is clear that $\mathcal{T}^\prime = \mathcal{T}\cap\Lambda^{u+1}_{S^\prime}\text{-mod},\mathcal{F}^\prime = \mathcal{F}\cap\Lambda^{u+1}_{S^\prime}\text{-mod}$. And $(\mathcal{T}^{\prime\prime},\mathcal{F}^{\prime\prime})$ has the similar property. So we only need to prove $\Delta^\prime_0\cup\Delta^\prime_1\cup\cdots\cup\Delta^\prime_u = \Delta^{\prime\prime}_0\cup\Delta^{\prime\prime}_1
\cup\cdots\cup\Delta^{\prime\prime}_v$.

For  convenience, we give the following notations for any given $i\geq0$:

$L^i_{S^\prime} = \langle Gen \mathrm{P}_{2j}(\Lambda^{2j}_{S^\prime})\mid 0\leq2j\leq max\{u,i\}\rangle $;

$R^i_{S^\prime} = \langle Cogen \mathrm{I}_{2j+1}(\Lambda^{2j+1}_{S^\prime})\mid 0\leq2j+1\leq max\{u,i\}\rangle $;

$L^i_{S^{\prime\prime}} = \langle Gen \mathrm{P}_{2j+1}(\Lambda^{2j+1}_{S^{\prime\prime}})\mid 0\leq2j+1\leq max\{v,i\}\rangle $;

$R^i_{S^{\prime\prime}} = \langle Cogen \mathrm{I}_{2j}(\Lambda^{2j}_{S^{\prime\prime}})\mid 0\leq2j\leq max\{v,i\}\rangle $.

We just prove $\Delta^\prime_0\cup\Delta^\prime_1\cup\cdots\cup\Delta^\prime_u\subseteq\Delta^{\prime\prime}_0\cup\Delta^{\prime\prime}_1
\cup\cdots\cup\Delta^{\prime\prime}_v$. For this, we just need to prove: $\forall i\geq0$, $L^{2i+1}_{S^\prime}\subseteq L^{2i+1}_{S^{\prime\prime}};R^{2i}_{S^\prime}\subseteq R^{2i}_{S^{\prime\prime}}$.

For $i = 0,R^0_{S^\prime} = \{0\}\subseteq Cogen \mathrm{I}_{0}(\Lambda^{0}_{S^{\prime\prime}}) = R^0_{S^{\prime\prime}}$, $L^1_{S^\prime} = Gen\mathrm{P}_0(\Lambda^0_{S^\prime})\subseteq Gen \mathrm{P}_1(\Lambda^1_{S^{\prime\prime}}) = L^1_{S^{\prime\prime}}$.

Now we assume the theorem holds for $i\leq k-1$. Then $\Lambda^{2k}_{S^{\prime\prime}}$ is a quotient algebra of $\Lambda^{2k-1}_{S^\prime}$ since $ \Delta^\prime_0\cup\Delta^\prime_1\cup\cdots\cup\Delta^\prime_{2k-2}\subseteq\Delta^{\prime\prime}_0\cup\Delta^{\prime\prime}_1
\cup\cdots\cup\Delta^{\prime\prime}_{2k-1}$. So $\Lambda^{2k}_{S^{\prime\prime}}\text{-mod}$ is a full subcategory of $\Lambda^{2k-1}_{S^\prime}\text{-mod}$.

Suppose $0\neq X\in \operatorname{add} \mathrm{I}_{2k-1}(\Lambda^{2k-1}_{S^\prime})$. So $X$ is $\Ext$-injective in $\Lambda^{2k}_{S^{\prime\prime}}\text{-mod}$. By  \text{torsion pair} $(\mathcal{F}\cap \Lambda^{2k}_{S^{\prime\prime}}\text{-mod},R^{2k-2}_{S^{\prime\prime}})$ on $\mathcal{F}$, there exists an exact sequence $0\rightarrow X_1\rightarrow X\rightarrow X_2\rightarrow 0$ such that $X_1\in \mathcal{F}\cap \Lambda^{2k}_{S^{\prime\prime}}\text{-mod}$ and $X_2\in R^{2k-2}_{S^{\prime\prime}}$. For every $Y\in \Lambda^{2k}_{S^{\prime\prime}}\text{-mod}$, applying Hom$_\Lambda(Y,-)$ to this exact sequence, we get an exact sequence:  $\Hom_\Lambda(Y,X_2) \rightarrow \Ext^1_\Lambda(Y,X_1) \rightarrow \Ext^1_\Lambda(Y,X)$. Since Ext$^1_\Lambda(Y,X) = 0$ and $Y\in \leftidx^\bot(R^{2k-2}_{S^{\prime\prime}})$, Ext$^1_\Lambda(Y,X_1) = 0$. So $X_1$ is $\Ext$-injective in $\Lambda^{2k}_{S^{\prime\prime}}\text{-mod}$. Thus $X_1\in \operatorname{add} \mathrm{I}_{2k}(\A^{2k}_{S^{\prime\prime}})$. So $X\in R^{2k}_{S^{\prime\prime}}$. Therefore,
$R^{2k}_{S^{\prime}}\subseteq R^{2k}_{S^{\prime\prime}}$, and similarly, we have $L^{2k+1}_{S^{\prime}}\subseteq L^{2k+1}_{S^{\prime\prime}}$.

\section{Examples}

\ \ \ \ \ In this section, we will use the results developed in the previous two sections to characterize \text{torsion pair}s on some particular module categories. Those results will be related to \cite{2}, \cite{3}, \cite{4}, \cite{5}, \cite{6}, \cite{11}. We always assume $K$ is a filed.  If Q is a quiver and $\Delta \in Q_0$ where $Q_0$ is the set of vertices of Q, then we denote the full sub-quiver of Q containing $\Delta$ by $Q(\Delta)$. We give the following definition.

\begin{defn}
Let Q be a quiver , $\{\Delta_0, \Delta_1, \dots, \Delta_m\}$  a tuple such that  $\Delta_i \subseteq Q_0, \ \Delta_i \bigcap \Delta_j  =  \phi \ \forall i \neq j, \Delta_0 \neq \phi$. If $\forall i >  0$ and $v \in \Delta_{2i + 1}$ there is a path from some vertex in $\Delta_{2i}$ to $v$ in the sub-quiver $Q(Q_0 - \Delta_0 - \Delta_1- \dots - \Delta_{2i - 1})$, and $\forall i >  0$ and $v \in \Delta_{2i}$ there is a path from $v$ to some vertex in $\Delta_{2i - 1}$in the sub-quiver $Q((Q_0 - \Delta_0 - \Delta_1- \dots - \Delta_{2i - 2})$. Then we call $\{\Delta_0, \Delta_1, \dots, \Delta_m\}$ is  a 2-type part partition of Q.  The following diagram shows the relation:
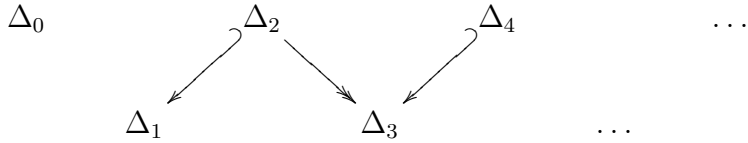
\[\xymatrix{ \Delta_0 &&\Delta_2 \ar@{->>}[rd] \ar@{_{(}->}[ld]&& \Delta_4 \ar@{_{(}->}[ld] && \dots\\
&\Delta_1 &&\Delta_3& &\dots }\]

 Dually, we we call $\{\Delta_0, \Delta_1, \dots, \Delta_m\}$ is a 2-type part partition of Q if $\forall i>  0$ and $v \in \Delta_{2i + 1}$ there is a path from $v$ to some vertex in $\Delta_{2i}$ in the sub-quiver $Q(Q_0 - \Delta_0 - \Delta_1- \dots - \Delta_{2i - 1})$, and $\forall i >  0$ and $v \in \Delta_{2i}$ there is a path from some vertex in $\Delta_{2i - 1}$ to $v$ in the sub-quiver $Q(Q_0 - \Delta_0 - \Delta_1- \dots - \Delta_{2i - 2})$.
The following diagram shows the relation:
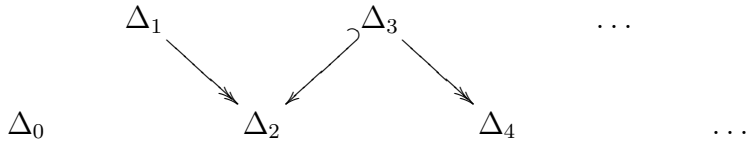
\[\xymatrix{ &\Delta_1 \ar@{->>}[rd]&&\Delta_3 \ar@{->>}[rd] \ar@{_{(}->}[ld]&& \dots\\
\Delta_0 &&\Delta_2 &&\Delta_4& &\dots  }\]

 Especially,  if $\forall i >  0, \Delta_{2i - 1}$ contains all sink points in $Q(Q_0 - \Delta_0 - \Delta_1- \dots - \Delta_{2i - 2})$, $\Delta
 _{2i}$ contains all source points in $Q(Q_0 - \Delta_0 - \Delta_1- \dots - \Delta_{2i - 1})$, then we call $\{\Delta_0, \Delta_1, \dots, \Delta_m\}$ is a strong 1-type part partition of Q. If $\forall i >  0, \Delta_{2i - 1}$ contains all source points in $Q(Q_0 - \Delta_0 - \Delta_1- \dots - \Delta_{2i - 2})$, $\Delta_{2i}$ contains all sink points in $Q(Q_0 - \Delta_0 - \Delta_1- \dots - \Delta_{2i - 1})$, then we call $\{\Delta_0, \Delta_1, \dots, \Delta_m\}$ is a strong 2-type part partition of Q.

 If $\Delta_0 \bigcup \Delta_1 \bigcup \dots \bigcup \Delta_m  =  Q_0$ we call $\{\Delta_0, \Delta_1, \dots, \Delta_m\}$ is a complete partition of Q.
\end{defn}

We have the following lemma.

\begin{lem}
Let Q be a acyclic quiver and $\{\Delta_0, \Delta_1, \dots, \Delta_m\}$ is a strong 1-type part partition of Q. Then $\{\Delta_0, \Delta_1, \dots, \Delta_m\}$ is a 1-type part partition of Q.

If $\{\Delta_0, \Delta_1, \dots, \Delta_m\}$ is a strong 2-type part partition of Q. Then $\{\Delta_0, \Delta_1, \dots, \Delta_m\}$ is a 2-type part partition of Q.
\end{lem}

For a quiver $Q$, we denote $\mathbf{E}(KQ)$ by  $\mathbf{E}(Q)$. Now we have the following theorem which is the path algebra's version of Theorem 3.11.

\begin{thm}
 Let Q be a acyclic quiver. Then we have a bijection between the set $(\mathcal{T}, \mathcal{F})$ which is a \text{torsion pair} on $KQ$-mod and the set of the pair $(\{\Delta_0, \Delta_1, \dots, \Delta_m\}$; $(\mathcal{T^\prime}, \mathcal{F^\prime}))$, where $\{\Delta_0, \Delta_1, \dots, \Delta_m\}$ is a  1-type part partition of Q and $(\mathcal{T^\prime}, \mathcal{F^\prime})\in \mathbf{E}(KQ(Q_0 - \Delta_0 - \Delta_1 - \dots - \Delta_m))$.
\end{thm}

The dual form of the theorem is similar, so we don't demonstrate here. Now let $A_n$ be the following quiver: $1 \ra 2 \ra 3 \ra \dots \ra n$. Applying the above theorem to the quiver $A_n$ , we have the following theorem.

\begin{thm} There exists a bijection between  \text{torsion pair}s on $KA_n$-mod and complete strong 1-type part partition sets of $A_n$.
\end{thm}

\noindent {\bf Proof:} It is easy to see $\mathbf{E}(KA_m)  =  \phi$ for every $m$. And a complete partition of Q is a 2-type part partition if and only if it is strong 1-type part partition. The rest is clear by the above theorem.\\

If we observe the bijection above, then we obtain some simple corollaries.
\begin{cor} Given a \text{torsion pair} $(\mathcal{T}, \mathcal{F})$ on $KA_n$-mod, then there exists a unique pair $(T, F)$ such that T, F are basic partial tilting modules,  $\#\operatorname{Ind}(T \bigoplus F)  =  n$, and $\mathcal{T}  =  Gen(T), \mathcal{F}  =  Cogen(F)$.
\end{cor}

\begin{cor}
If $\{\Delta_0, \Delta_1, \dots, \Delta_m\}$ is a complete strong 1-type part partition of $A_n$, then the corresponding \text{torsion pair} is induced by tilting modules if and only if $v_1 \in \Delta_0$.

If $\{\Delta_0, \Delta_1, \dots, \Delta_m\}$ is a complete strong 2-type part partition of $A_n$, then the corresponding \text{torsion pair} is induced by cotilting modules if and only if $v_n \in \Delta_0$.
\end{cor}

\begin{prop} The number of \text{torsion pair}s on $KA_n$ is the $(n + 1)-{th}$ Catalan number $C_{n + 1}  =  \frac{1}{n+2} {2n + 2 \choose n + 1}$.
\end{prop}

\noindent {\bf Proof:} Adding one vertex to $A_n$, then we have the quiver $A_{n + 1}: 1 \ra 2 \ra 3 \ra \dots \ra n \ra n + 1$. We have a \text{torsion pair} on $KA_{n + 1}$-mod: $(KA_n \text{-mod}, \mathcal{P}(KA_{n + 1}))$. So we have a bijection between \text{torsion pair}s on $KA_n \text{-mod}$ and \text{torsion pair}s induced by cotilting modules on $KA_{n + 1}$-mod by proposition 2.18. The number of \text{torsion pair}s induced by cotilting modules on $KA_{n + 1}$-mod is well known which is the $(n + 1)-{th}$ Catalan number(Lemma $A.1$ in \cite{11}).\\

\begin{defn}
 Suppose $\A$ is an artin algebra,   $\mathcal{C}$ is a full subcategory of $\A$-mod. If there exists a set of full subcategories $\{\mathcal{C}_i, i \in I\}$ of  $\mathcal{C}$ such that $\forall M \in \mathcal{C}$, there  uniquely exists a set of modules $M_{i_1} \in \mathcal{C}_1, M_{i_2} \in \mathcal{C}_2, \dots,  M_{i_n} \in \mathcal{C}_n$ where $i_1, i_2, \dots, i_n$ are mutually different such that $M \cong M_{i_1} \bigoplus M_{i_2} \bigoplus \dots \bigoplus M_{i_n}$, then we call $\mathcal{C}$ is the direct sum of $\{\mathcal{C}_i, i \in I\}$, and we denote $\mathcal{C}  =  \bigoplus_{i \in I} \mathcal{C}_i$.
 \end{defn}

 We have the following correspondence.

 \begin{lem}
 Suppose $\A$ is an artin algebra,   $\mathcal{C}$ is a full subcategory of $\A$-mod,  there exists a set of full subcategories $\{\mathcal{C}_i, i \in I\}$ of  $\mathcal{C}$ such that $\mathcal{C}  =  \bigoplus_{i \in I} \mathcal{C}_i$ and $\operatorname{Hom}(X, Y)  =  0$ for every $X \in \mathcal{C}_i, Y \in \mathcal{C}_j$ and $i \neq j$.
 Then there exists a bijection between \text{torsion pair}s on $\mathcal{C}$ and the tuple $\{(\mathcal{T}_i,   \mathcal{F}_i)\}_{i \in I}$ where $(\mathcal{T}_i,   \mathcal{F}_i)$ is a \text{torsion pair} on $\mathcal{C}_i$
  \end{lem}

\noindent {\bf Proof:}Given ${(\mathcal{T},   \mathcal{F})}$ a \text{torsion pair} on $\mathcal{C}$, then  ${(\mathcal{T} \bigcap \mathcal{C}_i,   \mathcal{F} \bigcap \mathcal{C}_i)}_{i \in I}$ is the corresponding tuple. Given  the tuple ${(\mathcal{T}_i,   \mathcal{F}_i)}_{i \in I}$ where $(\mathcal{T}_i,   \mathcal{F}_i)$,
then $(\bigoplus_{i \in I}\mathcal{T}_i, \bigoplus_{i \in I}\mathcal{F}_i)$ is the corresponding \text{torsion pair}.\\

Let $\tilde{A}_n$ be the following quiver with vertices $(\tilde{A}_n)_0  =  \{v_1, v_2, \dots v_n\}$:\\
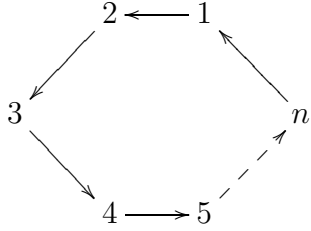
\[ \xymatrix{
  & 2  \ar[dl] &1\ar[l]\\
   3\ar[dr]&& &n\ar[ul]\\
   &4\ar[r] &5\ar@{-->}[ur]
}\]
Let $J$ be the ideal of $K\tilde{A}_n$ generated by all arrows. We call a finite-dimensional $K\tilde{A}_n$ module M is an ordinary module if there exists N such that $J^N M  =  0$. In this condition M is a $K\tilde{A}_n/J^N$ module. So if M is indecomposable, then it is uniserial and determined by its socle and length. Let $\mathcal{E}_n$ be the category of all ordinary modules. Then $\mathcal{E}_n$ is closed under submodules, quotients and extensions. We denote the simple module corresponding to the vertex $v_i$ by $S_i$.  We will give all \text{torsion pair}s on $\mathcal{E}_n$. For this we give the following definition which is introduced in \cite{2}.

\begin{defn}
Suppose $\Delta \in (\tilde{A}_n)_0$. let $Ray(\Delta)$ be the category of all modules with socle in $\operatorname{add} \bigoplus_{v_i \in \Delta} S_i$. let $Coray(\Delta)$ be the category of all modules with top in $\operatorname{add} \bigoplus_{v_i \in \Delta} S_i$.

For a subcategory $\mathcal{D}$ of $\mathcal{E}_n$. We denote $L_\mathcal{D}$ be the set of all vertices $v_i$ such that there are infinite indecomposable modules in $\mathcal{D}$ with $S_i$ as the top, $R_\mathcal{D}$ be the set of all vertices $v_j$ such that there are infinite indecomposable modules in $\mathcal{D}$ with $S_i$ as the socle  .
\end{defn}

By Definition 4.10 we have the following obvious lemma.
\begin{lem}
Suppose $\phi \neq \Delta \subseteq (\tilde{A}_n)_0$. Then $(Coray(\Delta), \tilde{A}_n((\tilde{A}_n)_0 - \Delta) \text{-mod})), (Q((\tilde{A}_n)_0 - \Delta) \text{-mod}, Ray(\Delta))$ are  two \text{torsion pair}s on $\mathcal{E}_n$.
\end{lem}

Now we give the following proposition.
\begin{prop}
Suppose $\phi \neq \Delta \subseteq (\tilde{A}_n)_0$. Then there is a bijection:\\
$\left( 1 \right)$  $\{(\mathcal{T^\prime},   \mathcal{F^\prime}): \text{ \text{torsion pair} on } \tilde{A}_n((\tilde{A}_n)_0 - \Delta)\text{-mod which is induced by cotilting modules} \} \\ \autorightleftharpoons{F}{$F^\prime$} \{(\mathcal{T},   \mathcal{F}): \text{ \text{torsion pair} on } \mathcal{E}_n\text{ such that } L_{\mathcal{T}}  =  \Delta \}$. In this condition $F((\mathcal{T^\prime},   \mathcal{F^\prime}))  =  (\langle Coray(\Delta), \mathcal{T^\prime}\rangle ,  \mathcal{F^\prime}), F^\prime((\mathcal{T},   \mathcal{F}))  =  (\mathcal{T} \bigcap \tilde{A}_n((\tilde{A}_n)_0 - \Delta)\text{-mod},   \mathcal{F})$
\\
$\left( 2 \right)$  $\{(\mathcal{T^\prime},   \mathcal{F^\prime}): \text{ \text{torsion pair} on } \tilde{A}_n((\tilde{A}_n)_0 - \Delta)\text{-mod which is induced by tilting modules} \} \\ \autorightleftharpoons{G}{$G^\prime$} \{(\mathcal{T},   \mathcal{F}): \text{ \text{torsion pair} on } \mathcal{C }\text{ such that }R_{\mathcal{F}} =  \Delta \}$. In this condition $G((\mathcal{T^\prime},   \mathcal{F^\prime}))  =  (\mathcal{T^\prime, \langle \mathcal{F^\prime}, \text{$Ray(\Delta)$}\rangle }), G^\prime((\mathcal{T},   \mathcal{F}))  =  (\mathcal{T},   \mathcal{F} \bigcap \tilde{A}_n((\tilde{A}_n)_0 - \Delta)\text{-mod})$.
\end{prop}

\noindent {\bf Proof:} We only proof (1) and (2) is similar.

$\left( 1 \right)$.  $\{(Coray(\Delta), \tilde{A}_n((\tilde{A}_n)_0 - \Delta) \text{-mod}), (\mathcal{E}_n, \{0\})\}$ is a 2-\text{torsion pair} seires on $\mathcal{E}_n$. Then   by Proposition 2.18, we have a bijection between \text{torsion pair} $(\mathcal{T},   \mathcal{F})$  on $\mathcal{E}_n$ such that  $Coray(\Delta) \subseteq \mathcal{T}$ and \text{torsion pair}s on  $\tilde{A}_n((\tilde{A}_n)_0 - \Delta) \text{-mod}$.  It is obvious in this condition $\Delta  =  L_\mathcal{T}$ if and only if in the corresponding \text{torsion pair} $(\mathcal{T^\prime},   \mathcal{F^\prime})$ on $\tilde{A}_n((\tilde{A}_n)_0 - \Delta)\text{-mod}$ $\mathcal{F^\prime}$ contains all projective modules which means it is induced by a cotilting module.\\

 The following lemma is from Corollary $4.5$ in \cite{2}.
\begin{lem}
Suppose $(\mathcal{T},   \mathcal{F}) \text{ is a \text{torsion pair}}$ on $\mathcal{E}_n$.  Then $L_{\mathcal{T}}, R_{\mathcal{F}}$ are not both empty.
\end{lem}

Now we have the following theorem which gives all \text{torsion pair}s on $\mathcal{E}_n$.

\begin{thm}
The following are all mutually different \text{torsion pair}s on $\mathcal{E}_n$ which are classified as two kinds. \\
$\left( 1 \right)$ $(Coray(\Delta) \bigoplus  \mathcal{T^\prime},  \mathcal{F^\prime})$ for some $\phi \neq \Delta \subseteq (\tilde{A}_n)_0$ and
$(\mathcal{T^\prime},   \mathcal{F^\prime})$ is a \text{torsion pair} on $ \tilde{A}_n((\tilde{A}_n)_0 - \Delta)\text{-mod}$ which is induced by cotilting modules.\\
$\left( 2 \right)$ $(\mathcal{T^\prime},  \mathcal{F^\prime} \bigoplus Ray(\Delta))$ for some $\phi \neq \Delta \subseteq (\tilde{A}_n)_0$ and
$(\mathcal{T^\prime},   \mathcal{F^\prime})$ is a \text{torsion pair} on $ \tilde{A}_n((\tilde{A}_n)_0 - \Delta)\text{-mod}$ which is induced by tilting modules.
\end{thm}

\noindent {\bf Proof:}
Suppose $(\mathcal{T},   \mathcal{F}) \text{ is a \text{torsion pair}}$ on $\mathcal{E}_n$ and $L_{\mathcal{T}} \neq \phi$. Then we know that \linebreak $Coray({L_{\mathcal{T}}) } \subseteq \mathcal{T}$ since $\mathcal{T}$ is closed under quotients. And for the first kind  it is obvious that $\langle Coray(\Delta), \mathcal{T^\prime}\rangle   =  Coray(\Delta) \bigoplus  \mathcal{T^\prime}$. The other is similar. \\

Since $\phi \neq \Delta$, we know $\tilde{A}_n((\tilde{A}_n)_0 - \Delta)\text{-mod}$ is a direct sum of module categories of $A_n$-type algebras. so by Lemma 4.9 the \text{torsion pair} is easily obtained. By the above theorem and the characterization of \text{torsion pair}s induced by tilting or cotilting modules on $A_n$-type algebras, we have the following bijection.

\begin{thm}
$\left(1\right)$ There is a  bijection between  the set of the \text{torsion pair}s $(\mathcal{T}, \mathcal{F})$ on $\mathcal{E}_n$ such that $L_\mathcal{T} \neq \phi$ and the set of the complete sets of $\tilde{A}_n$ $\{\Delta, \Delta_1, \dots, \Delta_m\}$ which is a strong 1-type part partition and $\Delta$ is not empty.\\
$\left(2\right)$ There is a  bijection between  the set of the \text{torsion pair}s $(\mathcal{T}, \mathcal{F})$ on $\mathcal{E}_n$ such that $R_\mathcal{F} \neq \phi$ and the set of the complete sest of $\tilde{A}_n$ $\{\Delta, \Delta_1, \dots, \Delta_m\}$ which is a strong 2-type part partition and $\Delta$ is not empty.
\end{thm}

\noindent {\bf Proof:} If $\{\Delta, \Delta_1, \dots, \Delta_m\}$ is a strong 1-type part partition, then $\{\Delta_1, \dots, \Delta_m\}$ is strong 1-type part partition in $\tilde{A}_n((\tilde{A}_n)_0 - \Delta)$. Then we get a \text{torsion pair}  $(\mathcal{T^\prime},  \mathcal{F^\prime})$ on $\tilde{A}_n((\tilde{A}_n)_0 - \Delta)$-mod which is induced by a cotilting module. Thus $(Coray(\Delta) \bigoplus \mathcal{T^\prime},  \mathcal{F^\prime})$ is the corresponding \text{torsion pair} on $\mathcal{E}_n$

If $\{\Delta, \Delta_1, \dots, \Delta_m\}$ is a strong 2-type part partition, then we get a \text{torsion pair} $(\mathcal{T^\prime},  \mathcal{F^\prime} \bigoplus Ray(\Delta))$ where $(\mathcal{T^\prime},  \mathcal{F^\prime})$ is a \text{torsion pair} on $\tilde{A}_n((\tilde{A}_n)_0 - \Delta)$-mod which is induced by a tilting module.

The rest is clear.\\

\section{\text{Torsion pair}s on hereditary algebras}

In this section we always assume $K$ is an algebraic closed field  and $Q$ is a acyclic quiver. We try to find a way  to obtain all \text{torsion pair}s on $KQ$-mod. This aim is also the motivation of the article. If $Q$ is not wild, we really get a way. If it is wild, the issue comes down to the \text{torsion pair}s on regular components of wild hereditary algebras. For this we denote the Auslander-Reiten translation by $\tau$, its quasi-inverse by $\tau^-$, the finite-dimensional projective $KQ$-module category by $\mathcal{P}(Q)$, the finite-dimensional injective $KQ$-module category by $\mathcal{I}(Q)$. The following two lemmas are well known.

\begin{lem}
Suppose $0 \ra A \ra B \ra C \ra 0$ is an exact sequence on
kQ-mod.Then\\
 $\left( 1 \right)$ If $\add A \bigcap \mathcal{P}(Q) = \{0\}$, then $0 \ra \tau A \ra \tau B \ra \tau C \ra 0$
is an exact sequence.\\
 $\left( 2 \right)$ If $\add C \bigcap \mathcal{I}(Q) = \{0\}$, then $0 \ra \tau^- A \ra \tau^- B \ra \tau^- C \ra 0$
is an exact sequence.\\
\end{lem}

\begin{lem}
Suppose $X, Y \in$ kQ-mod.\\
 $\left( 1 \right)$ If  $\add X \bigcap \mathcal{P}(Q) = \{0\}$, then $\Hom(X, Y) \cong \Hom(\tau X, \tau Y)$\\
 $\left( 2 \right)$ If $\add Y \bigcap \mathcal{I}(Q) = \{0\}$, then $\Hom(X, Y) \cong \Hom(\tau^- X, \tau^- Y)$
\end{lem}

We denote the set of \text{torsion pair}s on $KQ$-mod $(\mathcal{T}, \mathcal{F})$ such that $\mathcal{I}(Q) \subseteq \mathcal{T}$ by $\mathbf{F}_1(Q)$ and the set of \text{torsion pair}s on $KQ$-mod $(\mathcal{T}, \mathcal{F})$ such that $\mathcal{P}(Q) \subseteq \mathcal{F}$ by $\mathbf{F}_2(Q)$. And let $\mathbf{F}(Q) = \mathbf{F}_1(Q) \bigcup \mathbf{F}_2(Q)$. It is obvious that $\mathbf{E}(Q) = \mathbf{F}_1(Q) \bigcap \mathbf{F}_2(Q)$.  As a consequence of the above two lemmas, we have the following proposition. \\

\begin{prop}
Suppose there is no projective-injective $KQ$-module. Then
there is a one to one correspondence:\\
  $$\mathbf{F}_1(Q)  \autorightleftharpoons{$\sigma^-$}{$\sigma$} \mathbf{F}_2(Q)$$
such that
 $ \forall (\mathcal{T}^\prime, \mathcal{F}^\prime) \in \mathbf{F}_1(Q), \sigma^-(\mathcal{T}^\prime, \mathcal{F}^\prime) = (\tau^- \mathcal{T}^\prime,
\tau^- \mathcal{F}^\prime \bigoplus \mathcal{P}(Q));
  \forall (\mathcal{T}^{\prime\prime}, \mathcal{F}^{\prime\prime}) \in \mathbf{F}_2(Q), \linebreak  \sigma(\mathcal{T}^{\prime\prime}, \mathcal{F}^{\prime\prime}) = (\mathcal{I}(Q) \bigoplus \tau \mathcal{T}^{\prime\prime} ,
  \tau \mathcal{F}^{\prime\prime})$.
\end{prop}
\noindent {\bf Proof.}
We just prove that   $ \forall (\mathcal{T}^\prime, \mathcal{F}^\prime) \in \mathbf{F}_1, (\tau^- \mathcal{T}^\prime, \tau^- \mathcal{F}^\prime \bigoplus \mathcal{P}(Q))$ is a torsion pair on $KQ$-mod.

    By Lemma 5.2 (2), we know $\forall X \in \mathcal{T}^\prime, Y \in \mathcal{F}^\prime$, $\Hom(\tau^-X, \tau^-Y) \cong \Hom(X, Y) = \{0\}$. So the condition 1 in  the Definition 2.1 is satisfied.  By Lemma 5.1 (2), we know except projective modules, every indecomposable module has a suitable decomposition in $(\tau^- \mathcal{T}^\prime, \tau^- \mathcal{F}^\prime \bigoplus \mathcal{P}(Q))$. But for projective modules, the suitable decomposition is obvious. So the condition 2 in the Definition 2.1 is satisfied.\\

Just like the Auslander-Reiten translation, $\sigma^-$ and $\sigma$ also gives a translation on $F(Q)$. For every $(\mathcal{T}, \mathcal{F}) \in \mathcal{F}(Q), \text{if } \mathcal{I}(Q)
\subseteq \mathcal{T},$  then let $\sigma^-(\mathcal{T}, \mathcal{F}) = (\tau^- \mathcal{T}, \tau^- \mathcal{F} \bigoplus
\mathcal{P}(Q))$; if $\mathcal{P}(Q) \subseteq \mathcal{F}$, then let $\sigma(\mathcal{T}, \mathcal{F}) =
(\tau \mathcal{T} \bigoplus \mathcal{I}(Q), \tau \mathcal{F})$. The above proposition tells us that this translation defines
$\sigma$-obits for elements in $\mathbf{F}(Q)$. We use $[\mathcal{T},
\mathcal{F}]$ to denote the $\sigma$-obit of $(\mathcal{T}, \mathcal{F})$.

\begin{defn}
Suppose $(\mathcal{T}, \mathcal{F}) \in \mathbf{F}(Q)$. We call the elements in $[\mathcal{T}, \mathcal{F}] \bigcap (\mathbf{F}_2(Q) - \mathbf{F}_1(Q))$ source points of $[\mathcal{T}, \mathcal{F}]$, the elements in $[\mathcal{T}, \mathcal{F}] \bigcap (\mathbf{F}_1(Q) - \mathbf{F}_2(Q))$ sink points of $[\mathcal{T}, \mathcal{F}]$,  the elements in $[\mathcal{T}, \mathcal{F}] \bigcap \mathbf{F}_1(Q) \bigcap \mathbf{F}_2(Q)$ middle  points of $[\mathcal{T}, \mathcal{F}]$.
\end{defn}

The following corollary is obvious.
\begin{lem}
Suppose $(\mathcal{T}, \mathcal{F}) \in \mathbf{F}(Q)$.  Then $[\mathcal{T}, \mathcal{F}]$ has at most one source point and at most one sink  point. And
$[\mathcal{T}, \mathcal{F}] \bigcap \mathbf{F}_1(Q) \bigcap \mathbf{F}_2(Q) = [\mathcal{T}, \mathcal{F}] \bigcap \mathbf{E}(Q)$.
\end{lem}

We denote the preprojective component of $KQ$-mod by $\mathcal{P}_\infty(Q)$, the  preinjective component of $KQ$-mod by $\mathcal{I}_\infty(Q)$, the regular component of $KQ$-mod by $\mathcal{R}(Q)$.
\begin{thm}
Suppose $(\mathcal{T}, \mathcal{F}) \in \mathbf{F}(Q)$. Then\\
 $\left( 1 \right)$ $[\mathcal{T}, \mathcal{F}]$  has a source point but no sink point
 $\iff$ for every $(\mathcal{T}^\prime, \mathcal{F}^\prime) \in [\mathcal{T}, \mathcal{F}]$, $\mathcal{I}_\infty(Q) \bigcap \mathcal{F}^\prime \neq \phi$ and $\mathcal{P}_\infty(Q)
\subseteq \mathcal{F}^\prime$. \\
 $\left( 2 \right)$ $[\mathcal{T}, \mathcal{F}]$  has a sink point but no source point
$\iff$ for every $(\mathcal{T}^\prime, \mathcal{F}^\prime) \in [\mathcal{T}, \mathcal{F}]$, $\mathcal{P}_\infty(Q) \bigcap \mathcal{T} \neq \phi$ and
$\mathcal{I}_\infty(Q) \subseteq \mathcal{T}^\prime$. \\
 $\left( 3 \right)$ $[\mathcal{T}, \mathcal{F}]$  has a sink point and a source point $\iff$
for every $(\mathcal{T}^\prime, \mathcal{F}^\prime) \in [\mathcal{T}, \mathcal{F}]$,
$\mathcal{I}_\infty(Q) \bigcap \mathcal{F} \neq \phi$ and $\mathcal{P}_\infty(Q) \bigcap \mathcal{T}
\neq \phi$. \\
 $\left( 4 \right)$ $[\mathcal{T}, \mathcal{F}]$ has no sink point and no source point $\iff$ for every
$(\mathcal{T}^\prime, \mathcal{F}^\prime) \in [\mathcal{T}, \mathcal{F}]$, $\mathcal{I}_\infty(Q)
\subseteq \mathcal{T}^\prime$, and $\mathcal{P}_\infty(Q) \subseteq \mathcal{F}^\prime$.
\end{thm}

We denote the set of \text{torsion pair}s $(\mathcal{T}, \mathcal{F})$on $KQ$-mod such that $\mathcal{I}_\infty(Q)
\subseteq \mathcal{T}$, and $\mathcal{P}_\infty(Q) \subseteq \mathcal{F}$ by $\mathbf{H}(Q)$. So it is obvious that $\mathbf{H}(Q) \subseteq \mathbf{E}(Q)$.
We denote the set of \text{torsion pair}s on $\mathcal{R}(Q)$ by $\mathbf{R}(Q)$. We have the following obvious lemma.

\begin{lem}
There is a one to one correspondence:
$$\mathbf{H}(Q) \autorightleftharpoons{$F$}{$F^-$} \mathbf{R}(Q)$$
such that $\forall (\mathcal{T}, \mathcal{F}) \in \mathbf{H}(Q), F((\mathcal{T}, \mathcal{F})) = (\mathcal{T} \bigcap \mathcal{R}(Q), \mathcal{F} \bigcap \mathcal{R}(Q))$;  $\forall (\mathcal{T}^\prime, \mathcal{F})^\prime \in \mathbf{R}(Q), \linebreak F^-((\mathcal{T}^\prime, \mathcal{F}^\prime)) = (\mathcal{T}^\prime \bigoplus \mathcal{I}_\infty(Q), \mathcal{F}^\prime \bigoplus \mathcal{P}_\infty(Q)).$
\end{lem}

\begin{rem}
Suppose $(\mathcal{T}, \mathcal{F}) \in \mathbf{F}(Q)$ and  $[\mathcal{T}, \mathcal{F}]$ has at least  one  sink point or  one source point. We define the following operation $\Phi$:

Case $1$. If $[\mathcal{T}, \mathcal{F}]$ has a sink point, then we denote the sink point by $\Phi((\mathcal{T}, \mathcal{F}))$.

Case $2$.  If $[\mathcal{T}, \mathcal{F}]$ has a source point but no sink point, then we denote the source point by $\Phi((\mathcal{T}, \mathcal{F}))$.

For any torsion pair on $KQ$-mod we apply the operation in Theorem 3.11 and the operation $\Phi$ to it alternatively. At last  we get a new torsion pair on $KQ^\prime$-mod for some subquiver $Q^\prime$ of $Q$ such that the new torsion pair belongs to $\mathbf{H}(Q^\prime)$. This process is invertible by Theorem 3.11 and Proposition 5.3. So  by the above lemma if we  know  all torsion pairs on regular components for all subquivers, then we can construct all torsion pairs of $KQ$-mod.
\end{rem}

From now on we suppose $Q$ is a acyclic quiver with a Euclid ground graph. We start to find all the torsion pairs on $\mathbf{R}(Q).$ The following definition and two lemmas are from \cite{7}.

\begin{defn}
Suppose $X \in KQ$-mod . Then $Q$ is regular uniserial if there are regular submodules $0 = X_0
\subset X_1 \subset \dots \subset X_r = X$ and these are the only
regular submodules of X.
\end{defn}

\begin{lem}
If $\theta: X \ra Y$ with $X, Y$
regular $KQ$-modules, then $\Ima(\theta), \Ker(\theta)$ and $\Cok(\theta)$ are
regular.
\end{lem}

\begin{lem}
Every indecomposable regular $KQ$-module is regular universal.
\end{lem}

As an consequence we have
\begin{cor}
If $KQ$ is an Euclid-type algebra, X is a regular module, then the
quotient modules of X forms a chain: $X = X^r \twoheadrightarrow
\dots \twoheadrightarrow X^1 \twoheadrightarrow X^0$.
\end{cor}

\begin{cor}
Let $KQ$ be an Euclid-type algebra,  $f: X \ra Y$ is an injective morphism such that $X$ is a maximal regular submodule of the indecomposable regular module of $Y$. Then $f$ is an irreducible morphism.
\end{cor}

\noindent{\bf Proof.} $X$ is indecomposable by Lemma $5.11$. Suppose $\exists g:X \ra Z, h: Z \ra Z$ such that $f = hg$. Then by Lemma $5.11$, there is an indecomposable direct summand $Z^\prime$ such that $\exists g^\prime: X \ra Z^\prime, h: Z^\prime \ra Y$ such that $h^\prime g^\prime$ is an injective morphism. $Z^\prime$ is a regular module. So by Lemma $5.11, h^\prime$ is an injective morphism. Since  $X$ is a maximal regular submodule, $h^\prime$ is anisomorphism or  $g^\prime$ is an isomorphism.\\

Now Let $\mathcal{R}(Q) = \bigoplus_{i \in I}\mathcal{R}_i(Q)$
where $\{\mathcal{R}_i(Q), i \in I\}$ is the set of minimal additive categories containing  a connected  component in
AR-quiver of $KQ$. We denote the set of \text{torsion pairs} on  $\mathcal{R}_i(Q)$ by $\mathbf{R}_i(Q)$.  By Lemma 4.9, we have the following lemma.
\begin{cor}
There exists a bijection between $\mathbf{R}(Q)$ and the set of tuples $\{(\mathcal{T}_i,   \mathcal{F}_i)\}_{i \in I}$ with $(\mathcal{T}_i,   \mathcal{F}_i) \in \mathbf{R}_i(Q)$.
\end{cor}
\noindent {\bf Proof:} Let $X \in \mathcal{R}_i(Q)$. Then all  regular submodules and  all  regular quotient modules of $X$ are in $\mathcal{R}_i(Q)$ by the above corollary. So we know if $i \neq j$, then $\Hom(X, Y) = 0, \forall X \in \mathcal{R}_i(Q)$ and $Y \in \mathcal{R}_j(Q)$. The rest is clear by Lemma 4.9.\\

Now we start to demonstrate $\mathbf{R}_i(Q)$ . Suppose $\mathbf{R}_i(Q)$ has $n$ regular simple modules\cite{7}:
$S_1, S_2, \dots, S_{n - 1}$ where $S_{i + 1} = \tau S_i$. Let Let $\tilde{A}_n$ be the quiver in Section 4 and $S_1^\prime, S_2^\prime, \dots, S_n^\prime$ are the
correspondent simple modules to the vertices. Then we construct a
map: $\overline{F}(S_i^\prime) = S_i$. Then $\overline{F}$ induces a
one to one correspondence: $\mathcal{E}_n \ra \mathcal{R}_i(Q)$ such that
if $X$ $\in \mathcal{E}_n$ and is indecomposable with the length $m$ and top
$S_i^\prime$, then $F(X)$ is the indecomposable regular module with
the regular
length $m$ and top $S_i$. we have the following lemma.\\

 \begin{lem}
$\left( 1 \right)$ $\forall X, Y \in \mathcal{E}_n, \Hom(X, Y) = 0 \iff
\Hom(F(X), F(Y)) = 0$. \\
$\left( 2 \right)$ Suppose $Y \in \mathcal{E}_n$ and $X$ is a submodule of $Y$.
  Then $F(Y/X) = F(Y)/F(X)$.
 \end{lem}
\noindent {\bf Proof:} Clear by Lemma 5.11 and Corollary 5.12.

\begin{thm}
 F induces a one to one correspondent between the set of \text{torsion pair}s on
 $\mathcal{E}_n$ and  $\mathbf{R}_i(Q)$.
 \end{thm}
 \noindent {\bf Proof:} Clear by the above lemma.

\indent Fan Kong, Department of Mathematics, Shanghai Jiaotong  University, 200240 Shanghai, People's Republic of China.\\
\indent Email: Kongfan08@yahoo.com.cn\\

\indent Keyan Song, Department of Mathematics, Shanghai Jiaotong  University, 200240 Shanghai, People's Republic of China.\\
\indent Email: sky19840806@163.com\\

\indent Pu Zhang, Department of Mathematics, Shanghai Jiaotong  University, 200240 Shanghai, People's Republic of China.\\
\indent Email: pzhang@sjtu.edu.cn

\end{document}